\documentclass[review]{elsarticle}

\usepackage{lineno,hyperref}
\modulolinenumbers[5]

\journal{Journal of The Electrochemical Society}

\usepackage{mathrsfs}
\usepackage{amsfonts}
\usepackage{mathtools}
\usepackage{amsmath}
\usepackage{amsthm}
\usepackage{colortbl}
\usepackage{graphicx,verbatim}
\usepackage{subfig}

\newcommand{\re}{ {\rm Re} }
\newtheorem{thm}{Theorem}[section]

\newtheorem{lem}[thm]{Lemma}
\newtheorem{cor}[thm]{Corollary}
\newtheorem{asp}[thm]{Assumption}
\newtheorem{defn}[thm]{Definition}

\newcommand{\bd}[1]{\boldsymbol{#1}}

\newcommand{\pf}{\noindent {\bf Proof: }}
\newcommand{\eop}{{\hspace*{\fill}$\square$}}

\begin{document}

\begin{frontmatter}

\title{Efficient electrochemical model for lithium-ion cells}
\tnotetext[mytitlenote]{Fully documented templates are available in the elsarticle package on \href{http://www.ctan.org/tex-archive/macros/latex/contrib/elsarticle}{CTAN}.}

%\author{S. Afshar\fnref{myfootnote}}

\author[mymainaddress]{S. Afshar \corref{mycorrespondingauthor}}
\cortext[mycorrespondingauthor]{Corresponding author}
\ead{safshar@uwaterloo.ca}

\author[mymainaddress]{K. Morris}
\ead{kmorris@uwaterloo.ca}

\author[mysecondaryaddress]{A. Khajepour}
\ead{a.khajepour@uwaterloo.ca}

\address[mymainaddress]{Department of Applied Mathematics, University of Waterloo, Canada}
\address[mysecondaryaddress]{Department of Mechanical and Mechatronics Engineering, University of Waterloo, Canada}

\begin{abstract}
Lithium-ion batteries are used to  store  energy in electric vehicles. 
Physical models based on electro-chemistry  accurately predict  the cell dynamics, in particular the state of charge. However, these models are nonlinear partial differential equations coupled to algebraic equations, and they are computationally intensive.  Furthermore, a variable solid-state diffusivity model is recommended for cells with a lithium ion phosphate positive electrode to provide more accuracy. This variable structure adds more complexities to the model. 
However, a low-order model is required to represent the lithium-ion cells' dynamics for real-time applications.  
In this paper, a simplification of the electrochemical equations with variable solid-state diffusivity that preserves the key cells' dynamics is derived. The  simplified model is transformed into a numerically efficient fully dynamical form. It is proved that the simplified model is well-posed and can be approximated by a low-order finite-dimensional model. Simulations are very quick and show good agreement with experimental data. 
%Finally, a rate dependent diffusivity is  which improve the match to experimental data. 
 \end{abstract}


\end{frontmatter}

%\linenumbers

%================================================================
\section{Introduction and literature review}
%================================================================
%Cells are elementary units of  A battery pack, IS  the storage unit in electric/hybrid electric vehicles. 
Among different chemical compositions, a lithium-ion  chemistry is one of the most promising options for the batteries used for hybrid electric vehicles. 
High power and energy density, lack of memory effect, low self discharge, and high life cycle are some advantages of  lithium-ion chemistry in comparison to other cell chemistries \cite{est28-2013,est30-2012,est45-2012}. 
In particular, Lithium iron phosphate, LiFePO4 (LFP), offers the advantage of better lithium insertion over other alternatives. Its numerous features have drawn considerable interest. Some of these features are listed in \cite{2004-sirvi}.
%A battery management system is an important unit in the hybrid electric vehicles. The battery management system must ensure safety and efficient power utilization. 
Estimating state of charge (SOC), which determines the amount of deliverable energy, is critical for effective use of each cell and  for balancing the cells' state in a battery pack \cite{LIBA18-2008}. 
An accurate estimator that captures the cells' dynamics yet is simple enough for a real-time application is a important component of a battery management system.

Equivalent circuit models are frequently employed.  Simplicity and a relatively low number of parameters are the main advantages of these models \cite{mod5-2011}. Normally, the circuit includes a large capacitor or a voltage source to represent the open circuit potential (OCP) effect, and the rest of the circuit defines the cell's internal resistance and the effect of the cell's dynamics \cite{mod4-2011}. Different equivalent circuit models are introduced in \cite{est13-2004} and \cite{mod5-2011}. 
%Some examples of using equivalent circuit models in SOC estimation can be found in \cite{est17-2006,est25-2006,est12-2006,est64-2008,est20-2009,est42-2010,est19-2011,mod4-2011,est54,est41-2012}, and \cite{est45-2012}.

Electrochemical models, although more complex than the equivalent circuit models, 
have some advantages over other models in describing the cells' physical behavior. Including the effect of temperature and modeling the aging phenomenon, as well as other  inherent features of the lithium-ion batteries, is more feasible.   
%Unfortunately, the electrochemical models are not applicable in real-time applications because of their complexity.
 The electrochemical equations are nonlinear coupled partial differential equations. 
 These equations must be simplified without sacrificing their accuracy in order to obtain a model suitable for real-time applications.
Simplified and low-order models have been considered by many researchers. A review of most simplified electrochemical models is given in \cite{est29-2006} and \cite{LIBA24-2015}. 

A common simplification of the electrochemical equations is to assume that  there are only a finite number of particles along the electrodes. In a  single particle model each electrode is composed of a single spherical particle. In many cases this single-particle model provides good accuracy; see   \cite{est23-2006,est65-2010,est66-2010,August12-2,est69-2012,est68-2013, est28-2013,tanim2015aging} and \cite{tanim2015lithium}. 
In other situations a multiple particle model  with concentration-dependent solid diffusion coefficients
that considers the distribution of the particles in the electrodes provides better accuracy; see \cite{mp-r16}.
%This paper considers a model order reduction of this model such that the final model preserves most of the original equations' dynamics, and its outcome is close to experimental data of a commercial LFP cell. 

Several techniques have been developed to approximate the partial differential equations representing any simplified electrochemical model by ordinary differential equations. Laplace transforms and Pad\'e approximation are used in  \cite{LIBA12-2007,LIBA18-2008,LIBA25-2015}, and \cite{LIBA1-2011}. This approximation is also achieved via projection based techniques  such as proper orthogonal decomposition \cite{LIBA15-2008},  eigenfunctions of the solid diffusion equation  \cite{LIBA16-2006, LIBA2-2010}, and  orthogonal collocation \cite{Coordinate}. The approximation is derived using a polynomial approximation of the active material concentration in the solid phase in \cite{21-JES,Approximate}, and \cite{47_JES} and using Chebyshev polynomials in \cite{northrop2015}. A review of some approximation technique can be also found in \cite{33-JES} and \cite{northrop2014}.
These low-order models are introduced for a class of simplified models where the solid diffusion coefficient is often assumed to be constant. 

In practice, the solid phase diffusion coefficient is often a nonlinear function of active material concentration. A computationally efficient method, a control volume method, is developed in \cite{91234} for solving the diffusion equation with the variable diffusion equation. The approximation of the solid phase diffusion model with  the variable diffusion coefficient is also considered in \cite{urisanga2015} based on Lobatto IIIA quadrature to approximate the solid concentration. In this paper, eigenfunction based Galerkin collocation technique, which has shown an adequate result for constant diffusivity (\cite{33-JES}) and keeps key dynamical behaviour of the system, is extended to approximate the solid diffusion equation in which the diffusion coefficient is not constant.

Furthermore, the effect of porous electrode is often ignored in the electrochemical equations. Including the porous electrode model into the equations implies solving for a set of constraint equations simultaneously with the diffusion equations for the active material concentration. 
These constraint equations are coupled and nonlinear. 
A reduced order model is introduced in \cite{2002-Wang} in which a multi-scale model is developed to incorporate the pore level dynamics.
In most simplified models, the constraint equations are simplified by approximating the exchange current density with an average value; see  \cite{toward}, \cite{30-JES-2009}, \cite{LIBA20-2015, LIBA21-2015,LIBA22-2015,LIBA23-2015,LIBA25-2015}, and \cite{LIBA25-2015}. Linearizing the exchange current density term around some operating points is another method for approximation of the constraint equations (\cite{LIBA24a-2015} and \cite{LIBA1-2011}). However, in many applications of LiFePO4 (LFP)  cells, these approximations are not accurate; see, for instance, \cite{mpa-r17}. An efficient method for solving the constraint equations is considered here.

The main focus of this paper is developing a reduced order model for a full pseudo-two-dimensional electrochemical model with multiple variable solid-state diffusivity equations. A minimum number of approximations are introduced to facilitate the nonlinear analysis of the equations. These approximations are based on the physical properties of the system and have little effect on the results. First, it is proved that the solid and electrolyte potentials can be represented as differentiable functions of the solid and electrolyte concentrations as well as the input current. This representation is used to introduce a fully dynamical representation for the cell's dynamics. 
This simplified and transformed model were described in \cite{afshar2017fully}, where some simulation results were also provided. 

The cell's equations are also transformed into a state space form, which is proved to be well-posed. The state space representation is used to develop a low-order  model 
using the eigenfunction based Galerkin method that is shown to be efficient for real-time applications (\cite{33-JES}). It is proven that by adjusting the model order, an accuracy arbitrarily close to that of the original nonlinear partial differential equation model can be obtained. 

Finally, a fully dynamical low-order model is developed which is used in simulations. The simulation results show a good match to experimental data for different charging/discharging rates and profiles. It is also shown that the approximate solution converges as the order of approximation increases. Furthermore, the computation time  on a desktop computer  is  faster than the real time experimental time and comparable to the reported time for solving systems with constant diffusivity in the literature. Finally, it is described that the simulations match to the experimental data can be improved using a rate dependent diffusivity model. 

\section{Electrochemical model}

In this research, a lithium-ion cell with a positive electrode made of LFP material is considered. In LFP electrodes, the lithium insertion/deinsertion mechanism is a two phase process taking place between the lithium poor phase, $\text{Li}_{\epsilon}\text{FePO}_4$, and the lithium rich phase, $\text{Li}_{1-\epsilon}\text{FePO}_4$. The negative electrode is assumed to be a lithium foil.

A variable solid-state diffusivity model with a multiple particle size bins is used here. Details on this model can be found in \cite{mpa-r17}, and \cite{mp-r16}.
The battery cell's equations for the cell will be transformed to a state space representation. Let the number of particle size bins be $K=3$. Define  
\begin{align*}
&\mathcal{X}_1=\lbrace z\in\mathcal{L}^2([L_1,L]\times [0,R_1]):r_1^2z\in\mathcal{L}^2([L_1,L]\times [0,R_1]) \rbrace\\
&\mathcal{X}_2=\lbrace z\in\mathcal{L}^2([L_1,L]\times [0,R_2]):r_2^2z\in\mathcal{L}^2([L_1,L]\times [0,R_2]) \rbrace\\
&\mathcal{X}_3=\lbrace z\in\mathcal{L}^2([L_1,L]\times [0,R_3]):r_3^2z\in\mathcal{L}^2([L_1,L]\times [0,R_3]) \rbrace\\
&\mathcal{X}=\mathcal{L}^2(0,L)\times\mathcal{X}_1\times\mathcal{X}_2\times\mathcal{X}_3\\
&\mathcal{Y}=\mathcal{L}^2(0,L)\times\mathcal{L}^2(L_1,L)
\end{align*}
where $L_1=l_{sep}$, $L=l_{sep}+l_{cat}$.  The variable $c_e$ represents the electrolyte concentration, $c_{s,k}$ represent the solid concentration in each particle bin for $k=1,\ldots,3$, and $\varphi_e$ and $\varphi_s$ represent respectively the electrolyte and solid potential. 
Let $$\boldsymbol{c}=[c_1,c_2,c_3,c_4]^T=[c_e,c_{s,1},c_{s,2},c_{s,3}]^T\subseteq \mathcal{X}$$
 be the state vector and $\boldsymbol{\varphi}=[\varphi_1,\varphi_2]^T=[\varphi_e,\varphi_s]^T\subseteq \mathcal{Y}$ be the potential vector.
 
\begin{figure} [ht]
\centering
\includegraphics[width=.6\linewidth]{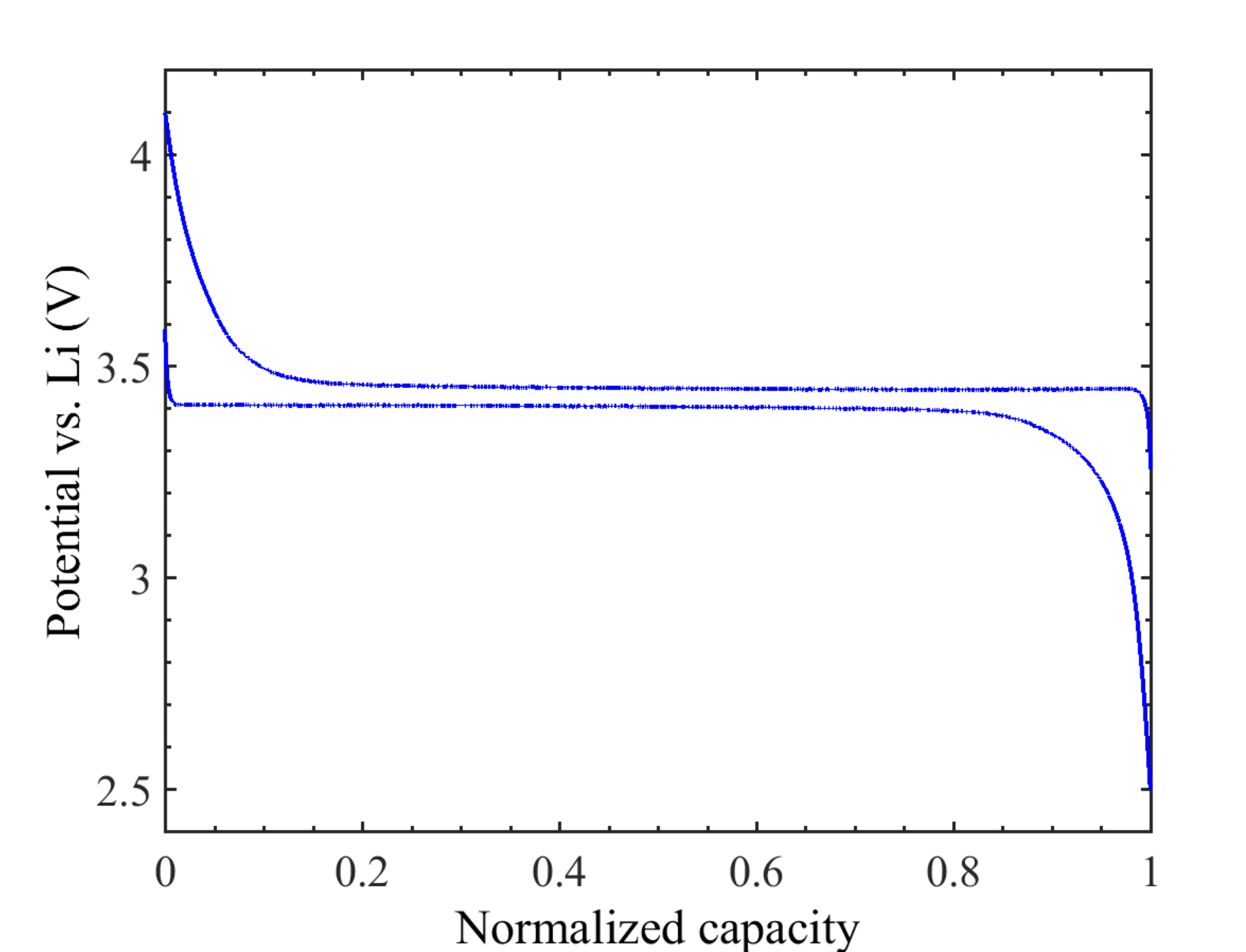}
\caption{OCP profile in a charging and discharging cycle (Laboratoire De R\'eactivit\'e Et Chimie Des Solides (LRCS) in Amiens, France).}
\label{ocp}
\end{figure}

Define 
\begin{equation}
\label{yk}
y_k=\text{sat}_y  \left(\frac{c_{k+1}\vert_{r_k=R_k}}{c_{max}} \right)
\end{equation}
where 
$$\text{sat}_y(s)=\frac{1}{1+\exp(-a_0s)}$$
for $a_0\in\mathbb{R}^+$ (see  Table \ref{T2}), 
The electrochemical reaction rate is defined as
\begin{align}
\label{ECD}
i_{k}(\boldsymbol{c},\boldsymbol{\varphi}) = \left\{ 
  \begin{array}{l l}
    2i_{0} \sinh\big(\frac{F \eta_k}{2RT}\big) & \quad 			\text{if $x\in [0,L_1]$}\\
    0 & \quad \text{if $x\in [L_1,L]$}\\
  \end{array} \right.
\end{align}
where
$$\eta_k=\varphi_2-\varphi_1-U(y_k)$$
and $U(\cdot)$ is the OCP term. Here
\begin{align*}
U(y_k)&=3.4510-0.009y_k
+0.6687\exp(-35y_k)-0.5\exp(-210(1-y_k))
\end{align*}
for the charging cycle, and 
\begin{align*}
U(y_k)&=3.4077-0.020269y_k
+0.5\exp(-200y_k)-0.9\exp(-30(1-y_k))
\end{align*}
for the discharging cycle.
The OCP profile has an important effect on the simulations and must be identified carefully. The OCP identification is based on the static performance and cannot be measured during the battery operation. Instead, an empirically derived relation are used. This empirical model is obtained through a curve fitting (the experimental data for OCP is shown in Figure \ref{ocp}; the source of the experimental data is quoted in the Simulation section).

The thermodynamic term or the activity correction factor $\alpha_k(\cdot)$ is defined in \cite{mp-r16} and was modified to be  a Fr\'echet differentiable function and fit experimental data as follows.
\begin{equation}
\label{acf}
\begin{aligned}
\alpha_k(c_{k+1})=
6\exp(-25y_k)+15\exp(-35(1-y_k))+0.3/(1+(y_k-0.5)^2).
\end{aligned}
\end{equation}
Note that
\begin{equation}
\label{bcon}
\delta_1\leq\vert\alpha_k(c_{k+1})\vert\leq\delta_2
\end{equation}
for $k=1,\ldots,3$ and $\delta_1,\delta_2\in\mathbb{R}^+$.

The cell governing equations are
\begin{align}
\label{E212a}
\frac{\partial}{\partial t}
\begin{bmatrix}
c_1\\
c_2\\
c_3\\
c_4\\
\end{bmatrix}
&=
\begin{bmatrix}
\frac{\partial}{\partial x}\big( D_{e}^{eff}\frac{\partial c_1}{\partial x}\big)+\frac{1-t^0_+}{F\epsilon}\sum a_k i_{k}(\boldsymbol{c},\boldsymbol{\varphi})\\
\frac{1}{r_1^2}\frac{\partial}{\partial r_1}\big( r_1^2\alpha_1(c_2)\mathcal{D}\frac{\partial c_2}{\partial r_1} \big)\\
\frac{1}{r_2^2}\frac{\partial}{\partial r_2}\big( r_2^2\alpha_2(c_3)\mathcal{D}\frac{\partial c_3}{\partial r_2} \big)\\
\frac{1}{r_3^2}\frac{\partial}{\partial r_3}\big( r_3^2\alpha_3(c_4)\mathcal{D}\frac{\partial c_4}{\partial r_3} \big)\\
\end{bmatrix}
\\[1ex]
\label{E212b}
0&=
\begin{bmatrix}
\frac{\partial}{\partial x}\big(k^{eff}\frac{\partial \varphi_1}{\partial x}\big)+
k^{eff}\frac{\partial}{\partial x}\big(\frac{2RT(1-t_+^0)}{Fc_1}\frac{\partial c_1}{\partial x}\big)+\sum a_k i_{k}(\boldsymbol{c},\boldsymbol{\varphi})\\
\frac{\partial}{\partial x}\big(\sigma^{eff}\frac{\partial \varphi_2}{\partial x}\big)
-\sum a_{k} i_{k}(\boldsymbol{c},\boldsymbol{\varphi})\\
\end{bmatrix}
\end{align}

The boundary conditions are
\begin{align}
\label{SB1p}
\frac{\partial c_1}{\partial x}\bigg\vert_{x=L}&=0\\
\label{SB2}
\frac{\partial c_{k+1}}{\partial r_k}\bigg\vert_{r_k=0}&=0,\quad k=1\dots 3\\
\label{SB5}
\varphi_{1}\bigg\vert_{x=0}&=0\\
\label{SB6}
\frac{\partial \varphi_{1}}{\partial x}\bigg\vert_{x=L}&=0\\
\label{SB4}
\frac{\partial\varphi_2}{\partial x}\bigg\vert_{x=L_1}&=0.
\end{align}
The controlled input is current $I(t)$,
\begin{align}
\label{SB1}
\epsilon_{sep}D_{e,sep}^{eff}\frac{\partial c_1}{\partial x}\bigg\vert_{x=0}&=-\frac{(1-t_+^0)I(t)}{F}\\
\label{SB7}
-\sigma^{eff}\frac{\partial \varphi_2}{\partial x}\bigg\vert_{x=L}&=I(t).
\end{align}
Also
\begin{equation}
\label{SB3}
\alpha_k(c_{k+1})\mathcal{D}\frac{\partial c_{k+1}}{\partial r_k}\bigg\vert_{r_k=R_k}=\frac{i_{k}}{F}(\boldsymbol{c},\boldsymbol{\varphi}),\quad k=1\dots 3.
\end{equation}

Finally, the solid potential in the negative electrode $\varphi_f$ satisfies 
$$I(t)=i_f\big(\frac{c_1}{c_{ini}}\big)^{1-\beta_f}\big(\exp\big(\frac{(1-\beta_f)F\varphi_f}{RT}\big)-\exp\big(\frac{\beta_f F\varphi_f}{RT}\big)\big)$$
where $c_{ini}$ is the initial value of the state variable $c_1$.

Some approximations are introduced to the model to facilitate nonlinear analysis of the equations including their well-posedness. First, 
the reaction rate is approximated; the variable $y_{k}$ defined in (\ref{yk}) is substituted by an average value.  Define
$$\bar{c}_{k+1}(x)=\int_0^{R_k}\delta(r_k-R_k)c_{k+1}(x,r_k)dr_k$$
and
$$\bar{y}_k=\text{sat}_y\big(\frac{\bar{c}_{k+1}}{c_{max}}\big)$$
where
\begin{align}
\label{del}
\delta(x-x_0) = \left\{ 
  \begin{array}{l l}
    \frac{1}{\epsilon_{0}x_0} & \quad 			\text{if $x\in [x_0-\epsilon_0x_0,x_0]$}\\
    0 & \quad \text{if $x\in [0,x_0-\epsilon_0x_0]$}\\
  \end{array} \right.
\end{align}
for some small $\epsilon_0>0$ (see Table \ref{T2}).
For parameters $b_0 , a_0$  (see  Table \ref{T2}), define  
$$\text{sat}(s)=\frac{2b_0}{1+\exp(-a_0s)}-b_0$$
and also define
$$\bar{\eta}_k=\varphi_2-\varphi_1-U(\bar{y}_k).$$
The exchange current density is approximated by
\begin{align}
\label{ECDA}
\bar{i}_k(\bd{c},\bd{\varphi}) = \left\{ 
  \begin{array}{l l}
    2i_{0} \sinh\big(\text{sat}\big(\frac{F \bar{\eta}_k}{2RT}\big)\big) & \quad 			\text{if $x\in [0,L_1]$}\\
    0 & \quad \text{if $x\in [L_1,L]$}\\
  \end{array} \right. \, . 
\end{align}
The argument of $\sinh(\cdot)$ is saturated in (\ref{ECDA}) to keep the electrochemical solution bounded. This constraint aligns with the physics of the system.

%\subsection{Numerically fast modeling}
A second  approximation is 
partially linearizing the constraint equations around the initial value of the electrolyte concentration $c_{ini}$.
The constraint equations become
\begin{equation}
\label{bat1}
0=
\begin{bmatrix}
\frac{\partial}{\partial x}\big(k^{eff}\frac{\partial \varphi_1}{\partial x}\big)+
k^{eff}\frac{\partial}{\partial x}\big(\frac{2RT(1-t_+^0)}{Fc_{ini}}\frac{\partial c_1}{\partial x}\big)+\sum a_k \bar{i}_{k}(\boldsymbol{c},\boldsymbol{\varphi})\\
\frac{\partial}{\partial x}\big(\sigma^{eff}\frac{\partial \varphi_2}{\partial x}\big)
-\sum a_{k} \bar{i}_{k}(\boldsymbol{c},\boldsymbol{\varphi})\\
\end{bmatrix}
\end{equation}
This approximation facilitates computation and also  guarantees that the system of constraint equation (\ref{E212b}) have a unique solution $\bd{\varphi}$ for every given state vector $\bd{c}.$
%The potential vector $\bd{\varphi}$ is a function of the state vector $\bd{c}. $
\begin{thm}
\label{IMP}
Define the operator $D\bd{\mathcal{O}}(\cdot):\mathcal{X}\times\mathcal{Y}\times\mathbb{R}^3\rightarrow\mathcal{Y}^{3\times 3}$ as
\begin{equation}
\label{bat4}
\begin{aligned}
&D\bd{\mathcal{O}}(\bd{c},\bd{\varphi},I(t),c_1(0),\varphi_2(L_1))=\\
&\begin{bmatrix}
\begin{matrix}
k^{eff}
+\sum_{k=1}^3a_k\int_0^x\int_0^y\frac{\partial \bar{i}_k(\bd{c}(s,r_k),\bd{\varphi}(s))}{\partial \varphi_1}dsdy\\
-x\int_{0}^{L}\sum_{k=1}^3a_k\frac{\partial \bar{i}_k(\bd{c}(s,r_k),\bd{\varphi}(s))}{\partial \varphi_1}ds
\end{matrix}
&
\begin{matrix}
\sum_{k=1}^3a_k\int_0^x\int_0^y\frac{\partial \bar{i}_k(\bd{c}(s,r_k),\bd{\varphi}(s))}{\partial \varphi_2}dsdy\\
-x\int_{0}^{L}\sum_{k=1}^3a_k\frac{\partial \bar{i}_k(\bd{c}(s,r_k),\bd{\varphi}(s))}{\partial \varphi_2}ds
\end{matrix}
&0\\
-\sum_{k=1}^3a_k\int_0^x\int_0^y\frac{\partial \bar{i}_k(\bd{c}(s,r_k),\bd{\varphi}(s))}{\partial \varphi_1}dsdy&\sigma^{eff}-\sum_{k=1}^3a_k\int_0^x\int_0^y\frac{\partial \bar{i}_k(\bd{c}(s,r_k),\bd{\varphi}(s))}{\partial \varphi_2}dsdy&\sigma^{eff}\\
\int_{0}^{L}\sum_{k=1}^3a_k\frac{\partial \bar{i}_k(\bd{c}(s,r_k),\bd{\varphi}(s))}{\partial \varphi_1}ds&
\int_{0}^{L}\sum_{k=1}^3a_k\frac{\partial \bar{i}_k(\bd{c}(s,r_k),\bd{\varphi}(s))}{\partial \varphi_2}ds&
0.
\end{bmatrix}
\end{aligned}
\end{equation}
If $D\bd{\mathcal{O}}(\cdot)$ is nonsingular at $[\bd{c}^*,\bd{\varphi}^*,I^*,c_0^*,\varphi_0^*]^T\in\mathcal{X}\times\mathcal{Y}\times\mathbb{R}^3$, the constraint equations (\ref{bat1}) have a unique solution such that the potential vector $\bd{\varphi}$ can be written as a Fr\'echet differentiable function of the state vector $\bd{c}$ and the input $I(t)$ in a neighborhood of this point. In other words, in some neighborhood of $[\bd{c}^*,\bd{\varphi}^*,I^*,c_0^*,\varphi_0^*]$,
\begin{equation}
\label{bat2}
\bd{\varphi}=\bd{\mathcal{R}}_{\bd{\varphi}}(\bd{c},c_1(0),I(t))
\end{equation} 
where $\bd{\mathcal{R}}_{\bd{\varphi}}(\cdot):\mathcal{X}\times\mathbb{R}^2\rightarrow\mathcal{Y}$ is a Fr\'echet differentiable function.
\end{thm}
\pf
In this proof, it is shown that $\bd{\varphi}$ is defined implicitly through the solution to an implicit equation $\bd{\mathcal{O}}(\bd{c},\bd{\varphi},I(t),c_1(0),\varphi_2(L_1))=\bd{0}$. It is proved that $\bd{\mathcal{O}}(\cdot)$ is Fr\'echet differentiable with derivative (\ref{bat4}). The proof is next a consequence of the Implicit Function Theorem \cite[Theorem 1.1.23]{IFT}. 

In the first step, define
$$\mathcal{O}_1(\cdot),\mathcal{O}_2(\cdot),\mathcal{O}_3(\cdot):\mathcal{X}\times\mathcal{Y}\times\mathbb{R}^3\rightarrow\mathcal{Y}$$
as
\begin{align*}
&\mathcal{O}_1(\bd{c},\bd{\varphi},I(t),c_1(0),\varphi_2(L_1))
=\frac{x}{k^{eff}}\int_{0}^{L}\sum_{k=1}^3a_k\bar{i}_{k}(\bd{c}(s,r_k),\bd{\varphi}(s))ds\\
&-\frac{2RT(1-t_+^0)}{Fc_{ini}}(c_1-c_1(0))-\frac{1}{k^{eff}}\int_0^x\int_0^y\sum_{k=1}^3a_k\bar{i}_{k}(\bd{c}(s,r_k),\bd{\varphi}(s))dsdy,\\
&\mathcal{O}_2(\bd{c},\bd{\varphi},I(t),c_1(0),\varphi_2(L_1))
=\varphi_2(L_1)+\frac{1}{\sigma^{eff}}\int_0^x\int_0^y\sum_{k=1}^3a_k\bar{i}_{k}(\bd{c}(s,r_k),\bd{\varphi}(s))dsdy,\\
&\mathcal{O}_3(\bd{c},\bd{\varphi},I(t),c_1(0),\varphi_2(L_1))
=I(t)+\int_{0}^{L}\sum_{k=1}^3a_k\bar{i}_{k}(\bd{c}(s,r_k),\bd{\varphi}(s))ds.
\end{align*}
Combined with the boundary conditions in (\ref{SB1p}), (\ref{SB5}), (\ref{SB6}), (\ref{SB4}), (\ref{SB1}), and (\ref{SB7}), the algebraic equation (\ref{bat1}) can be rewritten as
\begin{equation}
\label{bat3}
\begin{aligned}
&\varphi_1=\mathcal{O}_1(\bd{c},\bd{\varphi},I(t),c_1(0),\varphi_2(L_1))\\
&\varphi_2=\mathcal{O}_2(\bd{c},\bd{\varphi},I(t),c_1(0),\varphi_2(L_1))\\
&0=\mathcal{O}_3(\bd{c},\bd{\varphi},I(t),c_1(0),\varphi_2(L_1))
\end{aligned}
\end{equation}

Note that the functions $\bar{i}_k(\cdot)$ for $k=1,\ldots,3$ are Fr\'echet differentiable with respect to their arguments. This is due to the fact that $\text{sat}(\cdot)$, $\text{sat}_y(\cdot)$, and the empirical function chosen for OCP $U(\cdot)$, as well as the function $\sinh(\cdot)$ are Fr\'echet differentiable with respect to their arguments. Therefore, from the definition of $\bar{i}_k(\cdot)$ given by (\ref{ECDA}) and the chain rule, it can be concluded that $\bar{i}_k(\cdot)$ are Fr\'echet differentiable functions.  

Since integration is a linear operation, the fact that the functions $\bar{i}_k(\cdot)$ are Fr\'echet differentiable leads to the Fr\'echet differentiability of the functions $\mathcal{O}_1(\cdot)$, $\mathcal{O}_2(\cdot)$, and $\mathcal{O}_3(\cdot)$ with respect to $[c_2,\ldots,c_4]^T$ and $\bd{\varphi}$; these functions are linear and thus differentiable with respect to $(c_1-c_1(0))$, $\varphi_2(L_1)$, and $I(t)$. Define 
\begin{equation}
\label{O}
\begin{aligned}
\bd{\mathcal{O}}(\bd{c},\bd{\varphi}&,I(t),c_1(0),\varphi_2(L_1))=\\
&\begin{bmatrix}
k^{eff}(\varphi_1-\mathcal{O}_1(\bd{c},\bd{\varphi},I(t),c_1(0),\varphi_2(L_1)))\\
\sigma^{eff}(\varphi_2-\mathcal{O}_2(\bd{c},\bd{\varphi},I(t),c_1(0),\varphi_2(L_1)))\\
\mathcal{O}_3(\bd{c},\bd{\varphi},I(t),c_1(0),\varphi_2(L_1))
\end{bmatrix}
.
\end{aligned}
\end{equation}
The Fr\'echet derivative of the nonlinear operator $\bd{\mathcal{O}}(\cdot)$ (\ref{O}) with respect to the vector $[\bd{\varphi},\varphi_2(L_1)]^T$ is (\ref{bat4}). In addition, (\ref{bat3}) can be written as
$$\bd{\mathcal{O}}(\bd{c},\bd{\varphi},I(t),c_1(0),\varphi_2(L_1))=\bd{0}.$$
Now, by the Implicit Function Theorem and the assumption of $D\bd{\mathcal{O}}(\cdot)$ being nonsingular in some neighborhood of $[\bd{c}^*,\bd{\varphi}^*,I^*,c_0^*,\varphi_0^*]^T$, (\ref{bat2}) follows.
\eop

At this point, for the sake of simplicity and future use, $\bd{\mathcal{R}}_{\bd{\varphi}}(\cdot)$ in (\ref{bat2}) is approximated by 
\begin{equation}
\label{bat2ap}
\bd{\varphi}=\bar{\bd{\mathcal{R}}}_{\bd{\varphi}}(\bd{c},I(t))=\bd{\mathcal{R}}_{\bd{\varphi}}(\bd{c},\int_0^L\delta(x) c_1(x)dx,I(t))
\end{equation}
where
\begin{align}
\label{del0}
\delta(x) = \left\{ 
  \begin{array}{l l}
    \frac{1}{\epsilon_{0}L} & \quad 			\text{if $x\in [0,\epsilon_0L]$}\\
    0 & \quad \text{if $x\in [\epsilon_0L,L]$}\\
  \end{array} \right.
\end{align}
for some small $\epsilon_0\in\mathbb{R}^+$ given in Table \ref{T2}. This approximation is feasible due to the continuity of the electrolyte concentration.

Next, a new form of the constraint equations is achieved by taking the time differentiation of both sides of (\ref{bat2ap}). Along with (\ref{E212a}), differentiating (\ref{bat2ap}) results in
\begin{equation}
\label{bat5}
\frac{\partial\bd{\varphi}}{\partial t}=D\bar{\bd{\mathcal{R}}}_{\bd{\varphi}}(\bd{c},I(t))\big(\frac{\partial\bd{c}}{\partial t}\big)+\frac{\partial \bar{\bd{\mathcal{R}}}_{\bd{\varphi}}(\bd{c},I(t))}{\partial I}\frac{dI(t)}{dt}.
\end{equation}
The constraint equations (\ref{bat1}) are equivalent to the differential equations (\ref{bat5}). 
%With this setting, any standard technique of solving differential equations can be used without dealing with algebraic equations directly. 

Solving the differential equations (\ref{bat5}) requires the time derivative of  $I(t).$ This is accomplished  by using a saturated high-speed observer introduced in \cite{isi},
\begin{equation}
\label{bat6f}
\frac{d\hat{x}}{dt}=\bd{M}\hat{x}+\bd{L}I(t)
\end{equation}  
where $\hat{x}^T=[\hat{I},\hat{dI/dt}]$, and
$$
\bd{M}=
\begin{bmatrix}
-gh_1&1\\
-g^2h_0&0
\end{bmatrix},\quad
\bd{L}=
\begin{bmatrix}
-gh_1\\
-g^2h_0
\end{bmatrix}
$$
in which $g,h_0,h_1\in\mathbb{R}^+$ are tuning parameters. 

A third approximation of the cell's equations is made. 
Let both sides of (\ref{E212a}) followed by approximation (\ref{ECDA}) be multiplied by $\bd{w}=[w_1,\ldots,w_4]^T\in\mathcal{X}$ in the sense of the $\mathcal{X}$-inner product as follows:
\begin{equation}
\label{bat6}
\begin{aligned}
&\int_0^Lw_1(x)\frac{\partial c_1(x)}{\partial t}dx=\\
&\int_0^Lw_1(x)\big(\frac{\partial}{\partial x}\big( D_e^{eff}\frac{\partial  c_1(x)}{\partial x}\big)+\frac{1-t^0_+}{F\epsilon}\sum a_k \bar{i}_{k}(\boldsymbol{c}(x,r_k),\boldsymbol{\varphi}(x))\big)dx\\
&\int_{L_1}^L\int_0^{R_k}r_k^2w_{k+1}(x,r_k)\frac{\partial c_{k+1}(x,r_k)}{\partial t}dr_kdx=\\
&\int_{L_1}^L\int_0^{R_k}r_k^2w_{k+1}(x,r_k)
\frac{1}{r_k^2}\frac{\partial}{\partial r_k}\big( r_k^2\alpha_k(c_{k+1}(x,r_k))\mathcal{D}\frac{\partial c_{k+1}(x,r_k)}{\partial r_k} \big)dr_kdx\\
\end{aligned}
\end{equation}
for $k=1,\ldots,3.$
Now, applying integration by parts to (\ref{bat6}) and employing boundary conditions (\ref{SB1p}), (\ref{SB2}), (\ref{SB1}), and (\ref{SB3}) followed by approximation (\ref{ECDA}) lead to
\begin{equation}
\label{bat7}
\begin{aligned}
&\int_0^Lw_1(x)\frac{\partial c_1(x)}{\partial t}dx=\int_0^L\big(-\frac{\partial w_1(x)}{\partial x}\big( D_e^{eff}\frac{\partial c_1(x)}{\partial x}\big)\\
&+w_1(x)\frac{1-t^0_+}{F\epsilon}\sum a_k \bar{i}_{k}(\boldsymbol{c}(x,r_k),\boldsymbol{\varphi}(x))\big)dx
+\frac{1-t^0_+}{\epsilon F} w_1(0)I(t)\\
&\int_{L_1}^L\int_0^{R_k}r_k^2w_{k+1}(x,r_k)\frac{\partial c_{k+1}(x,r_k)}{\partial t}dr_kdx=-\int_{L_1}^L\int_0^{R_k}\frac{\partial w_{k+1}(x,r_k)}{\partial r_k}\\
&\big( r_k^2\alpha_k(c_{k+1}(x,r_k))\mathcal{D}\frac{\partial c_{k+1}(x,r_k)}{\partial r_k} \big)dr_kdx+\frac{R_k^2}{F}\int_0^Lw_{k+1}(R_k)\bar{i}_{k}(\boldsymbol{c}(x,r_k),\boldsymbol{\varphi}(x))dx.\\
\end{aligned}
\end{equation}

Next, (\ref{bat7}) is approximated by 
\begin{equation}
\label{bat8}
\begin{aligned}
&\int_0^Lw_1(x)\frac{\partial c_1(x)}{\partial t}dx=\int_0^L\big(-\frac{\partial w_1(x)}{\partial x}\big( D_e^{eff}\frac{\partial c_1(x)}{\partial x}\big)\\
&+w_1(x)\frac{1-t^0_+}{F\epsilon}\sum a_k \bar{i}_{k}(\bd{c}(x,r_k),\boldsymbol{\varphi}(x))\big)dx+\frac{1-t^0_+}{\epsilon F}\int_0^L\delta(x-L)w_1(x) I(t)dx\\
&\int_{L_1}^L\int_0^{R_k}r_k^2w_{k+1}(x,r_k)\frac{\partial c_{k+1}(x,r_k)}{\partial t}dr_kdx=\\
&-\int_{L_1}^L\int_0^{R_k}\frac{\partial w_{k+1}(x,r_k)}{\partial r_k}\big( r_k^2\alpha_k(c_{k+1}(x,r_k))\mathcal{D}\frac{\partial c_{k+1}(x,r_k)}{\partial r_k} \big)dr_kdx\\
&+\frac{R_k^2}{F}\int_0^L\int_0^{R_k}r_k^2\frac{\delta(r_k-R_k)}{r_k^2} w_{k+1}(x,r_k)\bar{i}_{k}(\bd{c}(x,r_k),\bd{\varphi}(x))dr_kdx.\\
\end{aligned}
\end{equation}
Using integration by parts in (\ref{bat8}), the battery equations can be transformed into 
\begin{equation}
\label{bal2-orig}
\frac{\partial\bd{c}}{\partial t}+\bd{\mathcal{A}}\bd{\mathcal{N}}(\bd{c})=\bd{\mathcal{R}}(\bd{c},t)+\bd{\mathcal{F}}\bd{u}(t)
\end{equation}
where 
%$\bd{\mathcal{A}}:\mathcal{X}\rightarrow\mathcal{X}$ is a linear operator, $\bd{\mathcal{R}}(.),\bd{\mathcal{N}}(.):\mathcal{X}\rightarrow\mathcal{X}$ are nonlinear Fr\'echet differentiable functions, $\bd{\mathcal{F}}:\mathbb{R}^q\rightarrow\mathcal{X}$ is a bounded linear input operator, and $\bd{u}(t)$ is the input vector defined as
\begin{equation}
\label{AOP}
\begin{aligned}
&\bd{\mathcal{A}}=+\lambda\bd{\mathcal{I}}\\
&-\begin{bmatrix}
\frac{\partial }{\partial x}\big( D_{e}^{eff}\frac{\partial }{\partial x}\big)&0&0&0\\
0&\frac{\mathcal{D}}{r_1^2}\frac{\partial }{\partial r_1}( r_1^2\frac{\partial }{\partial r_1} )&0&0\\
0&0&\frac{\mathcal{D}}{r_2^2}\frac{\partial }{\partial r_2}( r_2^2\frac{\partial }{\partial r_2} )&0\\
0&0&0&\frac{\mathcal{D}}{r_3^2}\frac{\partial }{\partial r_3}( r_3^2\frac{\partial }{\partial r_3} )\\
\end{bmatrix}
\end{aligned}
\end{equation}
with $\lambda>0$ is set such that $\bd{\mathcal{A}}$ is positive definite (this setting is required for future proofs), and
\begin{equation}
\label{dom}
\begin{aligned}
\mathcal{D}(\bd{\mathcal{A}})&=\lbrace \bd{c}\in\mathcal{X},[\frac{\partial c_1}{\partial x},\frac{\partial c_2}{\partial r_1},\ldots,\frac{\partial c_4}{\partial r_3}]^T\in\mathcal{X},\\
&[\frac{\partial^2 c_1}{\partial x^2},\frac{\partial^2 c_2}{\partial r_1^2},\ldots,\frac{\partial^2 c_4}{\partial r_3^2}]^T\in\mathcal{X}\\
&\frac{\partial c_1}{\partial x}(0)=\frac{\partial c_1}{\partial x}(L)=0,\, \mbox{and  }\\
& \frac{\partial c_{k+1}}{\partial r_{k}}(0)=\frac{\partial c_{k+1}}{\partial r_{k}}(R_k)=0 \,\mbox{for  }\ k=1\ldots 3\rbrace.
\end{aligned}
\end{equation}
\begin{align}
\label{ND}
&\bd{\mathcal{N}}(\bd{c})=
\begin{bmatrix}
c_1\\
\int_0^{r_1}\alpha_1(c_{2}(s,t))\frac{\partial c_{2}}{\partial r_1}(s,t)ds\\
\int_0^{r_2}\alpha_2(c_{3}(s,t))\frac{\partial c_{3}}{\partial r_2}(s,t)ds\\
\int_0^{r_3}\alpha_3(c_{4}(s,t))\frac{\partial c_{4}}{\partial r_3}(s,t)ds\\
\end{bmatrix}
\\
\label{RD}
&\bd{\mathcal{E}}(\bd{c},\bd{\varphi})=
\begin{bmatrix}
\frac{1-t^0_+}{F\epsilon}\sum a_k 
\bar{i}_{k}(\bd{c},\bd{\varphi})\\
\frac{R_1^2}{F}\frac{\delta(r_1-R_1)}{r_1^2} 
\bar{i}_{1}(\bd{c},\bd{\varphi})\\
\frac{R_2^2}{F}\frac{\delta(r_2-R_2)}{r_2^2} 
\bar{i}_{2}(\bd{c},\bd{\varphi})\\
\frac{R_3^2}{F}\frac{\delta(r_3-R_3)}{r_3^2} 
\bar{i}_{3}(\bd{c},\bd{\varphi})\\
\end{bmatrix}
+\lambda\bd{\mathcal{N}}(\bd{c})
\\
&\bd{\mathcal{B}}=
\begin{bmatrix}
\frac{1-t^0_+}{\epsilon F}\delta(x-L)\\
0\\
0\\
0\\
\end{bmatrix}\\
&\bd{\mathcal{R}}_{\bd{c}}(\bd{c},\bd{\varphi})=\bd{\mathcal{E}}(\bd{c},\bd{\varphi})-\bd{\mathcal{E}}(\bd{0},\bar{\bd{\mathcal{R}}}_{\bd{\varphi}}(\bd{0},I(t)))\\
&\bd{\mathcal{R}}(\bd{c},t)=\bd{\mathcal{R}}_{\bd{c}}(\bd{c},\bar{\bd{\mathcal{R}}}_{\bd{\varphi}}(\bd{c},I(t))),\\
&\bd{\mathcal{F}}=[\bd{\mathcal{B}},\bd{\mathcal{I}}],\quad \bd{u}(t)=[I(t),\bd{\mathcal{E}}(\bd{0},\bar{\bd{\mathcal{R}}}_{\bd{\varphi}}(\bd{0},I(t)))^T]^T
.
\end{align}

Thus, letting $D_{\bd{c}}$ represent  the Fr\'echet derivative with respect to $\bd{c}$, a fully dynamical representation equivalent to (\ref{bal2-orig}) is developed as
\begin{equation}
\label{CDNA-orig}
\begin{aligned}
\frac{\partial\bd{c}}{\partial t}&+\bd{\mathcal{A}}\bd{\mathcal{N}}(\bd{c})=\bd{\mathcal{R}}_{\bd{c}}(\bd{c},\bd{\varphi})+\bd{\mathcal{F}}\bd{u}(t)\\
\frac{\partial\bd{\varphi}}{\partial t}&=D_{\bd{c}}\bar{\bd{\mathcal{R}}}_{\bd{\varphi}}(\bd{c},I(t))\frac{\partial\bd{c}}{\partial t}+\frac{\partial \bar{\bd{\mathcal{R}}}_{\bd{\varphi}}(\bd{c},I(t))}{\partial I}\frac{dI(t)}{dt}\\
\end{aligned}
\end{equation}

The proof of the following lemma is a straightforward calculation.
\begin{lem}
\label{SADJ}
The linear operator $\bd{\mathcal{A}}:\mathcal{D}(\bd{\mathcal{A}})\in\mathcal{X}\rightarrow\mathcal{X}$ defined by (\ref{AOP}) is self-adjoint.  
\end{lem}

Furthermore, it can be easily checked that the inverse of the linear operator $\bd{\mathcal{A}}$ defined by (\ref{AOP}) is a double integral form with a bounded kernel; thus, it is a compact operator. This property along with the self-adjointness leads to the fact that the linear operator $\bd{\mathcal{A}}$ has eigenfunctions which are an orthogonal basis for the Hilbert space $\mathcal{X}$ \cite[theorem VIII.6]{robert}.

\section{Finite-dimensional approximation and well-posedness}
The electrochemical equations (\ref{bal2-orig}) is a special case of a general form 
\begin{equation}
\label{Sysg}
\frac{\partial\bd{z}}{\partial t}+\bd{\mathcal{A}}\bd{\mathcal{N}}(\bd{z})=\bd{\mathcal{R}}(\bd{z},t)+\bd{\mathcal{F}}\bd{u}(t)\\
\end{equation}
where $\bd{z}\in\mathcal{H}$ is the state vector, $\mathcal{H}$ is a Hilbert space and the state space, $\bd{\mathcal{R}}(\cdot):\mathcal{H}\rightarrow\mathcal{H}$ is a Fr\'{e}chet differentiable nonlinear operator with respect to $\bd{z}$ and strongly continuous with respect to $t$ that satisfies $\bd{\mathcal{R}}(\bd{0},t)=0$, and $\bd{\mathcal{N}}:\mathcal{H}\rightarrow\mathcal{H}$ is a Fr\'{e}chet differentiable nonlinear operator that satisfies $\bd{\mathcal{N}}(\bd{0})$. The operator $\bd{\mathcal{F}}$ is a bounded linear operator. 

The following assumptions are made for the general representation (\ref{Sysg}).
\begin{asp}
\label{asp1}
The control input $\bd{u}(t)$ is continuous in time and of bounded variation. In addition, there exist some $M_u\in\mathbb{R}^+$ such that $\Vert \bd{u}(t)\Vert\leq M_u$.
\end{asp}
\begin{asp}
\label{asp2}
The operator $\bd{\mathcal{A}}$ is assumed to be a self-adjoint closed operator with a compact inverse $\bd{\mathcal{A}}^{-1}$. It has also dense domain in Hilbert space $\overline{\mathcal{D}(\bd{\mathcal{A}})}=\mathcal{H}$ and is such that for some $\kappa>0$, $(\bd{\mathcal{A}}\bd{w},\bd{w})_{\mathcal{H}}\geq\kappa\Vert\bd{w}\Vert^2$ for every $\bd{w}\in\mathcal{D}(\bd{\mathcal{A}})$.
\end{asp}
\begin{asp}
\label{asp3}
The nonlinear operator $\bd{\mathcal{R}}(\cdot)$ is Lipschitz continuous on the Hilbert space $\mathcal{H}$. In other words, for every $\bd{w}_1,\bd{w}_2\in\mathcal{H}$, there exist a positive constant $L_R\in\mathbb{R}^+$ such that 
\begin{align*}
\Vert \bd{\mathcal{R}}(\bd{w}_1,t)-\bd{\mathcal{R}}(\bd{w}_2,t)\Vert_{\mathcal{H}}&\leq L_R \Vert \bd{w}_1-\bd{w}_2\Vert_{\mathcal{H}}.
\end{align*}
\end{asp} 
\begin{asp}
\label{gasp1}
The nonlinear operator $\bd{\mathcal{N}}(\cdot)$ is Fr\'echet differentiable and satisfies 
\begin{align*}
\epsilon_1\leq&\Vert D\bd{\mathcal{N}}(\bd{w}_2)\Vert\leq\epsilon_2 \\
\epsilon_1\Vert\bd{w}_1\Vert_{\mathcal{H}}^2\leq&(\bd{w}_1,D\bd{\mathcal{N}}(\bd{w}_2)\bd{w}_1)_{\mathcal{H}}=(D\bd{\mathcal{N}}(\bd{w}_2)\bd{w}_1,\bd{w}_1)_{\mathcal{H}}\leq \epsilon_2\Vert\bd{w}_1\Vert_{\mathcal{H}}^2
\end{align*}
for every $\bd{w}_1,\bd{w}_2\in\mathcal{H}$ and some $\epsilon_1,\epsilon_2>0$.
\end{asp}
\begin{asp}
\label{gasp2}
The linear operator $\bd{\mathcal{A}}$ and the nonlinear operator $\bd{\mathcal{N}}(\cdot)$ satisfy 
$$(\bd{w},\bd{\mathcal{A}}\bd{\mathcal{N}}(\bd{w}))_{\mathcal{H}}=(\bd{\mathcal{A}}\bd{\mathcal{N}}(\bd{w}),\bd{w})_{\mathcal{H}}\geq\epsilon_3\Vert\bd{w}\Vert^2$$
for every $\bd{w}\in\mathcal{D}(\bd{\mathcal{A}}^{1/2})$ such tat $\bd{\mathcal{N}}(\bd{w})\in\mathcal{D}(\bd{\mathcal{A}})$ and some $\epsilon_3>0$.
\end{asp}

The linear operator $\bd{\mathcal{A}}$ can also be used to define a new Hilbert space with more smoothness properties. Before the normed space of interest can be defined, the concept of evolution triple and duality pairing are introduced first. This definition will be used in next section to prove the well-posedness of the equations.
\begin{defn}
\label{Gelf}
(Duality Pairing, \cite[Definition 3.4.3]{gelfan})\\
Let $\mathcal{V}\subseteq\mathcal{H}$ be a linear space whose dual space is denoted by $\mathcal{V}'$. The triple $(\mathcal{V},\mathcal{H},\mathcal{V}')$ is called an evolution triple if it satisfies the following conditions:
\begin{itemize}
\item the linear space $\mathcal{V}$ is a separable and reflexive Banach space.
\item the linear space $\mathcal{H}$ is a separable Hilbert space.
\item for $\mathcal{V}\subseteq\mathcal{H}\subseteq\mathcal{V}'$, $\mathcal{V}$ is dense and continuously embedded in $\mathcal{H}$.
\end{itemize}
The duality pairing between $\mathcal{V}$ and $\mathcal{V}'$ is denoted by $\langle.,.\rangle_{\mathcal{V}',\mathcal{V}}$ and defined as a continuous extension of the inner product on the Hilbert space $\mathcal{H}$, denoted by $(\cdot,\cdot)_{\mathcal{H}}$.
\end{defn}

From Assumption \ref{asp2}, $\bd{\mathcal{A}}^{1/2}$ is a well-defined positive definite operator; thus it is possible to define a Hilbert space $\mathcal{V}=\mathcal{D}(\bd{\mathcal{A}}^{1/2})$ with norm $\Vert \bd{\mathcal{A}}^{1/2} \cdot\Vert_{\mathcal{H}}$. With this setting, $\mathcal{V}$ is dense in the Hilbert space $\mathcal{H}$ and $\bd{\mathcal{A}}^{1/2}$ defines an isomorphism between $\mathcal{V}$ and $\mathcal{H}$ since it is a bounded linear operator from $\mathcal{V}$ to $\mathcal{H}$ with bounded linear inverse from $\mathcal{H}$ to $\mathcal{V}$. Therefore, $(\mathcal{V},\mathcal{H},\mathcal{V}')$ is a evolution triple and a duality pairing can be defined as in Definition \ref{Gelf}. 

Furthermore, for every $w\in\mathcal{V}$, 
$$(\bd{\mathcal{A}}^{1/2}\cdot,w)_{\mathcal{H}}:\mathcal{V}\rightarrow\mathbb{C}$$
is a linear functional with domain $\mathcal{V}$ dense in $\mathcal{H}$; thus, it can be extended uniquely to the Hilbert space $\mathcal{H}$ by Hahn-Banach theorem. This extension is the dual pairing between $\mathcal{V}$ and $\mathcal{V}'$. 
Respectively, from the definition of duality pairing, Definition \ref{Gelf}, for $\bd{w}_1\in\mathcal{H}$ and $\bd{w}_2,\bd{w}_3\in\mathcal{V}$,
\begin{equation}
\label{eq1}
\begin{aligned}
(\bd{w}_1,\bd{w}_2)_{\mathcal{H}}&=\langle\bd{w}_1,\bd{w}_2\rangle_{\mathcal{V}',\mathcal{V}} \, ,\\
(\bd{\mathcal{A}}^{1/2}\bd{w}_2,\bd{\mathcal{A}}^{1/2}\bd{w}_3)_{\mathcal{H}}&=\langle\bd{\mathcal{A}}\bd{w}_2,\bd{w}_3\rangle_{\mathcal{V}' \, , \mathcal{V}}.
\end{aligned}
\end{equation}

\begin{defn}
\label{SSol}
(Strong solution, \cite{Temam})\\ 
A strong solution to (\ref{Sysg}) is an element $\bd{z}\in\mathcal{H}$ which
\begin{itemize}
\item is strongly continuous and differentiable in time for almost every $t\in[0,t_f]$ with respect to $\mathcal{H}$-norm topology,
\item satisfies $\bd{z}(0)=\bd{z}_0$ for the initial condition $\bd{z}_0\in\mathcal{H}$,
\item and satisfies equation (\ref{Sysg}) for almost every $t\in[0,t_f]$.
\end{itemize} 
\end{defn}

Given assumption \ref{asp2}, the eigenfunctions of the linear operator $\bd{\mathcal{A}}$ provide an orthogonal basis for the Hilbert space $\mathcal{H}$ \cite[theorem VIII.6]{robert}.
The eigenfunctions $\bd{v}_i$ of the linear operator $\bd{\mathcal{A}}$ and the Galerkin method are used to define a finite-dimensional Hilbert space $\mathcal{H}_N$. An orthonormal projection onto the finite-dimensional Hilbert space $\mathcal{H}_N$ is defined by
$$\bd{\mathcal{P}}_N\bd{z}=\sum_{i=1}^N z_i\bd{v}_i.$$

Let the system's state be approximated by $\bd{z}_N=\bd{\mathcal{P}}_N\bd{z}$. The reduced order system is defined as 
\begin{equation}
\label{res1}
\frac{\partial\bd{z}_N}{\partial t}=\bd{\mathcal{A}}\bd{\mathcal{N}}_N(\bd{z}_N)+\bd{\mathcal{R}}_N(\bd{z}_N,t)+\bd{\mathcal{F}}_N\bd{u}(t)\\
\end{equation}
where 
\begin{align*}
\bd{\mathcal{N}}_N(\cdot)&=\bd{\mathcal{P}}_N\bd{\mathcal{N}}(\cdot)\\
\bd{\mathcal{A}}_N&=\bd{\mathcal{P}}_N\bd{\mathcal{A}}\\
\bd{\mathcal{R}}_N(\cdot)&=\bd{\mathcal{P}}_N\bd{\mathcal{R}}(\cdot)\\
\bd{\mathcal{F}}_N&=\bd{\mathcal{P}}_N\bd{\mathcal{F}}.\\
\end{align*}

%-------------------------------------------------------------------------------------------------
The following Lemma shows the boundedness of the solution to (\ref{res1}).
\begin{lem}
\label{bdd}
Let the system (\ref{Sysg}) satisfy Assumption \ref{asp1}-\ref{asp3}, \ref{gasp1}, and \ref{gasp2}. Suppose that $\bd{\mathcal{N}}(\bd{z}(x,0))\in\mathcal{V}$. The solution to (\ref{res1}) on every bounded time interval $[0,t_f]$ is bounded;
\begin{align}
\Vert\bd{z}_N(t)\Vert_{\mathcal{H}}&\leq M_{c,0}\\
\label{bdd1}
\Vert\bd{\mathcal{A}}^{1/2}\bd{\mathcal{N}}_N(\bd{z}_N(t))\Vert_{\mathcal{H}}&\leq M_{c,1}\\
\label{bdd2}
\int_0^{t_f}\Vert\bd{\mathcal{A}}\bd{\mathcal{N}}_N(\bd{z}_N(t))\Vert^2_{\mathcal{H}}dt&\leq M_{c,2}\\
\end{align} 
for $M_{c,0},M_{c,1},M_{c,2}\in\mathbb{R}^+$ independent of $N$.
\end{lem}
\pf
First,
from Assumption \ref{gasp1} and Mean Value Theorem \cite{Mean}[Theorem 7.6-1], it is concluded that $\bd{\mathcal{N}}(\cdot)$ is Lipschitz continuous. In other words, for every $\bd{w}_1,\bd{w}_2\in\mathcal{H}$ and some $L_{\mathcal{N}}>0$,
\begin{equation}
\label{LPZ}
\begin{aligned}
\Vert\bd{\mathcal{N}}(\bd{w}_2)-\bd{\mathcal{N}}(\bd{w}_1)\Vert_{\mathcal{H}}\leq &L_{\mathcal{N}}\Vert\bd{w}_2-\bd{w}_1\Vert_{\mathcal{H}}\\
\end{aligned}
\end{equation}
Note that
\begin{equation}
\label{bdd3}
\bd{\mathcal{R}}(\bd{0},t)=\bd{0},\quad \bd{\mathcal{N}}(\bd{0})=\bd{0}.
\end{equation}
Furthermore, by Assumption \ref{gasp2},
\begin{equation}
\label{bdd4}
\begin{aligned}
(\bd{z}_N,\bd{\mathcal{A}}\bd{\mathcal{N}}_N(\bd{z}_N))_{\mathcal{H}}=(\bd{\mathcal{A}}\bd{\mathcal{N}}(\bd{z}_N),\bd{z}_N)_{\mathcal{H}}\geq 0.
\end{aligned}
\end{equation}
Let both sides of (\ref{res1}) be multiplied by $\bd{w}\in\mathcal{H}$; 
\begin{equation}
\label{bat8a}
\begin{aligned}
\big(\bd{w},\frac{\partial\bd{z}_N}{\partial t}\big)_{\mathcal{H}}+(\bd{w},\bd{\mathcal{A}}\bd{\mathcal{N}}_N(\bd{z}_N))_{\mathcal{H}}=(\bd{w},\bd{\mathcal{R}}_N(\bd{z}_N,t)+\bd{\mathcal{F}}_N\bd{u}(t))_{\mathcal{H}}.\\
\end{aligned}
\end{equation}
Similarly,
\begin{equation}
\label{bat8ab}
\begin{aligned}
\big(\frac{\partial\bd{z}_N}{\partial t},\bd{w}\big)_{\mathcal{H}}+(\bd{\mathcal{A}}\bd{\mathcal{N}}_N(\bd{z}_N),\bd{w})_{\mathcal{H}}=(\bd{\mathcal{R}}_N(\bd{z}_N,t)+\bd{\mathcal{F}}_N\bd{u}(t),\bd{w})_{\mathcal{H}}.\\
\end{aligned}
\end{equation}
Next, replacing $\bd{w}$ by $\bd{z}_N$ in (\ref{bat8a}) and (\ref{bat8ab}) and adding the resulting equations yield
\begin{equation}
\label{bat8ac}
\begin{aligned}
\frac{d\Vert\bd{z}_N\Vert^2_{\mathcal{H}}}{dt}+2(\bd{\mathcal{A}}\bd{\mathcal{N}}_N(\bd{z}_N),\bd{z}_N)_{\mathcal{H}}=2\re(\bd{\mathcal{R}}_N(\bd{z}_N,t)+\bd{\mathcal{F}}_N\bd{u}(t),\bd{z}_N)_{\mathcal{H}}.\\
\end{aligned}
\end{equation}
Employing (\ref{bdd4}), the Lipschitz continuity (\ref{LPZ}), and (\ref{bdd3}) as well as using  Cauchy Schwarz and Young's inequality in (\ref{bat8ac}) leads to
\begin{equation}
\label{bdd5}
\frac{d\Vert\bd{z}_N(t)\Vert^2_{\mathcal{H}}}{dt}\leq L_1\Vert\bd{z}_N(t)\Vert^2_{\mathcal{H}}+L_2
\end{equation}
where
$$L_1=2L_R+1,\quad L_2=\Vert\bd{\mathcal{F}}\Vert M_u$$
and $M_u$ is the upper bound of $\bd{u}(t)$. Integrating inequality (\ref{bdd5}) results in
\begin{equation}
\label{bdd6}
\Vert\bd{z}_N(t)\Vert^2_{\mathcal{H}}\leq\Vert\bd{z}_N(0)\Vert^2_{\mathcal{H}}\exp(L_1t)+\frac{L_2(\exp(L_1t)-1)}{L_1}\leq M_{c,0}
\end{equation}
for some $M_{c,0}>0$.

Now, let both sides of (\ref{res1}) be first operated by $D\bd{\mathcal{N}}_N(\bd{z}_N)$, the Fr\'echet derivative of $\bd{\mathcal{N}}_N(\cdot)$, and then multiplied by $\bd{\mathcal{A}}\bd{\mathcal{N}}_N(\bd{z}_N)$ in the sense of the Hilbert space inner product; it is derived from following the same procedure as before that
\begin{equation}
\label{bdd7}
\begin{aligned}
\big(\bd{\mathcal{A}}\bd{\mathcal{N}}_N(\bd{z}_N),D\bd{\mathcal{N}}_N(\bd{z}_N)\frac{\partial\bd{z}_N}{\partial t}\big)_{\mathcal{H}}&+\big(D\bd{\mathcal{N}}_N(\bd{z}_N)\frac{\partial\bd{z}_N}{\partial t},\bd{\mathcal{A}}\bd{\mathcal{N}}_N(\bd{z}_N)\big)_{\mathcal{H}}=\\
&-2(\bd{\mathcal{A}}\bd{\mathcal{N}}_N(\bd{z}_N),D\bd{\mathcal{N}}_N(\bd{z}_N)\bd{\mathcal{A}}\bd{\mathcal{N}}_N(\bd{z}_N))_{\mathcal{H}}+\\
&2\re(\bd{\mathcal{A}}\bd{\mathcal{N}}_N(\bd{z}_N),D\bd{\mathcal{N}}_N(\bd{z}_N)(\bd{\mathcal{R}}_N(\bd{z}_N,t)+\bd{\mathcal{F}}_N\bd{u}(t)))_{\mathcal{H}}.
\end{aligned}
\end{equation}
Note that from Fr\'echet differentiability of $\bd{\mathcal{N}}(\cdot)$, for $\bd{h}\in\mathcal{H}$,
\begin{align*}
\frac{\Vert\bd{\mathcal{P}}_N(\bd{\mathcal{N}}(\bd{z}_N+\bd{h})
-\bd{\mathcal{N}}(\bd{z}_N)-D\bd{\mathcal{N}}(\bd{z}_N)\bd{h})\Vert_{\mathcal{H}}}{\Vert\bd{h}\Vert_{\mathcal{H}}}\rightarrow 0
\end{align*}
when $\Vert\bd{h}\Vert_{\mathcal{H}}\rightarrow 0$; therefore,
\begin{equation}
\label{DN1}
D\bd{\mathcal{N}}_N(\bd{z}_N)=\bd{\mathcal{P}}_ND\bd{\mathcal{N}}(\bd{z}_N).
\end{equation}
From (\ref{DN1}) and the fact that $\bd{\mathcal{A}}\bd{\mathcal{N}}_N(\bd{z}_N)\in\mathcal{H}_N$,
it is concluded that
\begin{align*}
(\bd{\mathcal{A}}\bd{\mathcal{N}}_N(\bd{z}_N),D\bd{\mathcal{N}}_N(\bd{z}_N)\bd{\mathcal{A}}\bd{\mathcal{N}}_N(\bd{z}_N))_{\mathcal{H}}=
(\bd{\mathcal{A}}\bd{\mathcal{N}}_N(\bd{z}_N),D\bd{\mathcal{N}}(\bd{z}_N)\bd{\mathcal{A}}\bd{\mathcal{N}}_N(\bd{z}_N))_{\mathcal{H}},
\end{align*}
and, from Assumption \ref{gasp1}
\begin{equation}
\label{DN2a}
\begin{aligned}
-(\bd{\mathcal{A}}\bd{\mathcal{N}}_N(\bd{z}_N),D\bd{\mathcal{N}}_N(\bd{z}_N)\bd{\mathcal{A}}\bd{\mathcal{N}}_N(\bd{z}_N))_{\mathcal{H}}
\leq-\epsilon_1 \Vert\bd{\mathcal{A}}\bd{\mathcal{N}}_N(\bd{z}_N)\Vert^2_{\mathcal{H}}.
\end{aligned}
\end{equation}
Similarly,
\begin{equation}
\label{DN2b}
\begin{aligned}
-(D\bd{\mathcal{N}}_N(\bd{z}_N)\bd{\mathcal{A}}\bd{\mathcal{N}}_N(\bd{z}_N),\bd{\mathcal{A}}\bd{\mathcal{N}}_N(\bd{z}_N))_{\mathcal{H}}
\leq-\epsilon_1 \Vert\bd{\mathcal{A}}\bd{\mathcal{N}}_N(\bd{z}_N)\Vert^2_{\mathcal{H}}.
\end{aligned}
\end{equation}
Substituting (\ref{DN2a}) and (\ref{DN2b}) into (\ref{bdd7}) and employing Cauchy Schwarz inequality; Young's inequality; and Assumption \ref{asp1}, \ref{asp3}, and \ref{gasp1} in (\ref{bdd7}) lead to 
\begin{equation}
\label{bdd9}
\begin{aligned}
\frac{d\Vert\bd{\mathcal{A}}^{1/2}\bd{\mathcal{N}}_N(\bd{z}_N(t))\Vert^2_{\mathcal{H}}}{dt}&
\leq-L_3\Vert\bd{\mathcal{A}}\bd{\mathcal{N}}_N(\bd{z}_N(t))\Vert^2_{\mathcal{H}}+L_4
\end{aligned}
\end{equation}
where 
$$L_3=2\epsilon_1-\beta_5,\quad L_4=\frac{1}{\beta_5}(\beta_3L_RM_{c,0}+\Vert\bd{\mathcal{F}}\Vert M_u)^2$$
and $\beta_5$, which comes from Young's inequality, is set such that $L_3>0$. 
Since $$-L_3\Vert\bd{\mathcal{A}}\bd{\mathcal{N}}_N(\bd{z}_N)\Vert^2_{\mathcal{H}}<0,$$ by integrating (\ref{bdd9}) and employing (\ref{bdd6}) on the bounded time interval $[0,t_f]$ the second boundedness result is achieved as
\begin{equation}
\label{bdd10}
\Vert\bd{\mathcal{A}}^{1/2}\bd{\mathcal{N}}_N(\bd{z}_N(t))\Vert^2_{\mathcal{H}}\leq\Vert\bd{\mathcal{A}}^{1/2}\bd{\mathcal{N}}_N(\bd{z}_N(0))\Vert^2_{\mathcal{H}}+L_4t_f\leq M_{c,1}
\end{equation}
for $M_{c,1}\in\mathbb{R}^+$.

Integrating (\ref{bdd9}) and considering the boundedness given by (\ref{bdd10}) lead to
\begin{equation}
\label{bdd11}
\begin{aligned}
\int_0^{t_f}\Vert\bd{\mathcal{A}}\bd{\mathcal{N}}_N(\bd{z}_N(t))\Vert^2_{\mathcal{H}}dt&\leq 
\frac{L_4}{L_3}t_f\\
&+\frac{1}{L_3}(\Vert\bd{\mathcal{A}}^{1/2}\bd{\mathcal{N}}_N(\bd{z}_N(0))\Vert^2_{\mathcal{H}}-\Vert\bd{\mathcal{A}}^{1/2}\bd{\mathcal{N}}_N(\bd{z}_N(t_f))\Vert^2_{\mathcal{H}})\\
&\leq M_{c,2}
\end{aligned}
\end{equation}
for some $M_{c,2}\in\mathbb{R}^+$.
\eop
%-------------------------------------------------------------------------------------------------

\begin{thm}
\label{ABDD}
Let the assumptions of Lemma \ref{bdd} be satisfied. The system (\ref{Sysg}) has at least one strong solution $\bd{z}\in\mathcal{L}^{2}([0,t_f];\mathcal{V})\cap\mathcal{L}^{\infty}([0,t_f];\mathcal{H})$. Furthermore, the approximation error $\bd{e}_N=\bd{z}-\bd{z}_N$ is bounded, and the sequence $\bd{e}_N$ admits a subsequence converging to zero in $\mathcal{L}^2([0,t_f];\mathcal{H})$ as $N$ goes to infinity. 
\end{thm}
\pf
It can be concluded from Lemma \ref{bdd} that the sequence $\bd{z}_N$ stays in a bounded set in  $\mathcal{L}^{\infty}([0,t_f];\mathcal{V})$ and thus in $\mathcal{L}^{2}([0,t_f];\mathcal{V})\cap\mathcal{L}^{\infty}([0,t_f];\mathcal{H})$. It is also concluded that $\bd{\mathcal{N}}_N(\bd{z}_N)$ stays in a bounded set in $\mathcal{L}^{2}([0,t_f];\mathcal{V})$. By Banach-Alaoglu theorem \cite{sell}, there exists a subsequence $\bd{z}_M$ and $\bd{\mathcal{N}}_M(\bd{z}_M)$ such that
\begin{equation}
\label{st1a}
\begin{aligned}
&\bd{z}_M(t)\rightarrow\bd{z}^*(t)\,\mbox{weakly in }\,\mathcal{L}^{2}([0,t_f];\mathcal{V})\\
&\bd{z}_M(t)\rightarrow\bd{z}^*(t)\,\mbox{in weak-star topology in }\,\mathcal{L}^{\infty}([0,t_f];\mathcal{H}),
\end{aligned}
\end{equation} 
and 
\begin{equation}
\label{st1aa}
\bd{\mathcal{N}}_M(\bd{z}_M(t))\rightarrow\bd{w}^*(t)\,\mbox{weakly in }\,\mathcal{L}^{2}([0,t_f];\mathcal{V})
\end{equation}
for $\bd{z}^*(t)\in\mathcal{L}^{2}([0,t_f];\mathcal{V})\cap\mathcal{L}^{\infty}([0,t_f];\mathcal{H})$ and $\bd{w}^*(t)\in\mathcal{L}^{2}([0,t_f];\mathcal{V})$ since $\mathcal{L}^2([0,t_f];\mathcal{V})$ and $\mathcal{L}^2([0,t_f];\mathcal{H})$ are complete with respect to weak topology.  
From (\ref{res1}), Lipschitz continuity (\ref{LPZ}), boundedness of $\bd{u}(t)$, and Lemma \ref{bdd}, it is concluded that the sequence $d\bd{z}_M(t)/dt$ stays in a bounded set in $\mathcal{L}^{2}([0,t_f];\mathcal{H})$. Therefore, by \cite[Theorem III.2.1]{TemamR},
\begin{equation}
\label{st3}
\bd{z}_M(t)\rightarrow\bd{z}^*(t)\,\mbox{strongly in }\,\mathcal{L}^{2}([0,t_f];\mathcal{H}).
\end{equation} 

Note that from (\ref{st1aa}), it can be concluded that for $\bd{w}(t)\in\mathcal{L}^{2}([0,t_f];\mathcal{D}(\bd{\mathcal{A}}))$,
\begin{align*}
\int_0^{t_f}(\bd{\mathcal{A}}^{\frac{1}{2}}\bd{\mathcal{N}}_M(\bd{z}_M(t)),\bd{\mathcal{A}}^{\frac{1}{2}}\bd{w}(t))_{\mathcal{H}}dt\rightarrow
\int_0^{t_f}(\bd{\mathcal{A}}^{\frac{1}{2}}\bd{w}^*(t),\bd{\mathcal{A}}^{\frac{1}{2}}\bd{w}(t))_{\mathcal{H}}dt;
\end{align*}
thus,
\begin{equation}
\label{wst1}
\int_0^{t_f}(\bd{\mathcal{N}}_M(\bd{z}_M(t)),\bd{\mathcal{A}}\bd{w}(t))_{\mathcal{H}}dt\rightarrow
\int_0^{t_f}(\bd{w}^*(t),\bd{\mathcal{A}}\bd{w}(t))_{\mathcal{H}}dt.
\end{equation}
In addition, by (\ref{st3}),
\begin{equation}
\label{wst2}
\int_0^{t_f}(\bd{\mathcal{N}}_M(\bd{z}_M(t)),\bd{\mathcal{A}}\bd{w}(t))_{\mathcal{H}}dt\rightarrow
\int_0^{t_f}(\bd{\mathcal{N}}(\bd{z}^*(t)),\bd{\mathcal{A}}\bd{w}(t))_{\mathcal{H}}dt.
\end{equation}
Since $\bd{\mathcal{A}}$ has a bounded linear inverse by Assumption \ref{asp2}, it is onto $\mathcal{H}$. Therefore, the convergence results (\ref{wst1}) and (\ref{wst2}) are satisfied for every $\bar{\bd{w}}(t)=\bd{\mathcal{A}}\bd{w}(t)\in\mathcal{L}^{\infty}([0,t_f];\mathcal{H})$. Therefore, by uniqueness of the limit in weak topology, $\bd{w}^*(t)=\bd{\mathcal{N}}(\bd{z}^*(t))$, and
\begin{equation}
\label{wst4}
\bd{\mathcal{N}}_M(\bd{z}_M(t))\rightarrow\bd{\mathcal{N}}(\bd{z}^*(t))\,\mbox{weakly in }\,\mathcal{L}^{2}([0,t_f];\mathcal{V}).
\end{equation}

Now, multiplying both sides of (\ref{bat8ab}) by a smooth function $\phi(t)$ with $\phi(t_f)=0$, employing (\ref{eq1}), and integrating by part
the resulting equation with respect to time yield
\begin{equation}
\label{bat8b}
\begin{aligned}
-\int_0^{t_f}\big((\bd{z}_M(t),\bd{w})_{\mathcal{H}}\frac{d\phi(t)}{d t}&+(\bd{\mathcal{A}}^{\frac{1}{2}}\bd{\mathcal{N}}_M(\bd{z}_M(t)),\bd{\mathcal{A}}^{\frac{1}{2}}\bd{w})_{\mathcal{H}}\phi(t)\big)dt\\
&=\int_0^{t_f}(\bd{\mathcal{R}}_M(\bd{z}_M(t),t)+\bd{\mathcal{F}}_M\bd{u}(t),\bd{w})_{\mathcal{H}}\phi(t)dt\\
&+(\bd{z}_M(0),\bd{w})_{\mathcal{H}}\phi(0).\\
\end{aligned}
\end{equation}
For $\bd{w}\in\mathcal{D}(\bd{\mathcal{A}}^{1/2})$, passing the limits (\ref{st1a}), (\ref{st3}), (\ref{wst4}), and the limit
$$\bd{z}_M(0)\rightarrow\bd{z}(0)\,\mbox{strongly in }\,\mathcal{H}$$
to (\ref{bat8b}) and using Assumption \ref{asp3} lead to
\begin{equation}
\label{bat8bb}
\begin{aligned}
-\int_0^{t_f}\big((\bd{z}^*(t),\bd{w})_{\mathcal{H}}\frac{d\phi(t)}{d t}&+(\bd{\mathcal{A}}^{\frac{1}{2}}\bd{\mathcal{N}}(\bd{z}^*(t)),\bd{\mathcal{A}}^{\frac{1}{2}}\bd{w})_{\mathcal{H}}\phi(t)\big)dt\\
&=\int_0^{t_f}(\bd{\mathcal{R}}(\bd{z}^*(t),t)+\bd{\mathcal{F}}\bd{u}(t),\bd{w})_{\mathcal{H}}\phi(t)dt+
(\bd{z}(0),\bd{w})_{\mathcal{H}}\phi(0).\\
\end{aligned}
\end{equation}
Finally, integrating (\ref{bat8bb}) by parts results in
\begin{equation}
\label{bat8bbb}
\begin{aligned}
\int_0^{t_f}\frac{d}{d t}(\bd{z}^*(t),\bd{w})_{\mathcal{H}}\phi(t)dt=&-\int_0^{t_f}(\bd{\mathcal{A}}^{\frac{1}{2}}\bd{\mathcal{N}}(\bd{z}^*(t)),\bd{\mathcal{A}}^{\frac{1}{2}}\bd{w})_{\mathcal{H}}\phi(t)dt\\
&+\int_0^{t_f}(\bd{\mathcal{R}}(\bd{z}^*(t),t)+\bd{\mathcal{F}}\bd{u}(t),\bd{w})_{\mathcal{H}}\phi(t)dt.\\
\end{aligned}
\end{equation}
Using (\ref{eq1}) in (\ref{bat8bbb}) yields to
\begin{equation}
\label{abdd1}
\frac{d}{dt}\langle\bd{z}^*,\bd{w}\rangle_{\mathcal{V}',\mathcal{V}}=\langle-\bd{\mathcal{A}}\bd{\mathcal{N}}(\bd{z}^*)+\bd{\mathcal{R}}(\bd{z}^*,t)+\bd{\mathcal{F}}\bd{u}(t),\bd{w}\rangle_{\mathcal{V}',\mathcal{V}}
\end{equation}
which is valid in distribution sense on $[0,t_f]$.
Since 
\begin{align*}
\bd{z}^*(t)&\in\mathcal{L}^2([0,t_f];\mathcal{H}),\\
-\bd{\mathcal{A}}\bd{\mathcal{N}}(\bd{z}^*(t))+\bd{\mathcal{R}}(\bd{z}^*(t),t)+\bd{\mathcal{F}}\bd{u}(t)&\in\mathcal{L}^2([0,t_f];\mathcal{H}),
\end{align*}
by \cite[Lemma II.3.1]{Temam} and from (\ref{abdd1})
$$\frac{\partial\bd{z}^*(t)}{\partial t}\in\mathcal{L}^2([0,t_f];\mathcal{H})$$
and $\bd{z}^*(t)$ satisfies (\ref{Sysg}) almost every where. Furthermore, by \cite[Lemma II.3.1]{Temam}, $\bd{z}^*$ equals almost every where to a continuous function from $[0,t_f]$ to $\mathcal{H}$; thus, it is a strong solution to (\ref{Sysg}) by Definition \ref{SSol}.
\eop

Now, the electrochemical equations are shown to satisfy  Assumptions \ref{asp2}-\ref{gasp2}.
\begin{cor}
\label{LM}
In the system (\ref{bal2-orig}), let the input signal $\bd{u}$ satisfies Assumption \ref{asp1}, and $\bd{\mathcal{N}}(\bd{c}(0))\in\mathcal{D}(\bd{\mathcal{A}}^{1/2})$ holds for the initial condition where $\bd{\mathcal{A}}$ and $\bd{\mathcal{N}}(\cdot)$ are defined respectively by (\ref{AOP}) and (\ref{ND}). Then, the system (\ref{bal2-orig}) has a strong solution $\bd{c}$. Furthermore, for the state vector $\bd{z}=\bd{c}$, the system can be approximated by finite-dimensional equations with the same form as (\ref{res1}) whose solutions $\bd{c}_N$ admit a convergent subsequence in $\mathcal{L}^2([0,t_f];\mathcal{X})$ where $[0,t_f]$ is a finite time interval.
\end{cor}
\pf
First, by Lemma \ref{SADJ}, $\bd{\mathcal{A}}$ defined by (\ref{AOP}) is self-adjoint.
Furthermore, as mentioned before, it can be easily checked that the inverse of the linear operator $\bd{\mathcal{A}}$  is a double integral form with a bounded kernel; thus, it is a compact operator. This property along with the self-adjointness leads to the fact that the linear operator $\bd{\mathcal{A}}$ satisfies Assumption \ref{asp2} \cite[theorem VIII.6]{robert}. 

%-------------------------------------------------------------------------------------------------
Next, it is proved that the nonlinear operators $\bd{\mathcal{N}}(\cdot)$ and $\bd{\mathcal{R}}(\cdot)$ satisfy Assumptions \ref{asp3},\ref{gasp1}, and \ref{gasp2}. First, it can be concluded from the definition of $\bd{\mathcal{N}}(\cdot)$ and chain rule theorem \cite[Theorem 3.2.1]{Chain} that 
\begin{equation}
\label{NN}
D\bd{\mathcal{N}}(\bd{c})=
\begin{bmatrix}
\mathcal{I}&0&0&0\\
0&\alpha_1(c_2)&0&0\\
0&0&\alpha_2(c_3)&0\\
0&0&0&\alpha_3(c_4)\\
\end{bmatrix}
.
\end{equation}
Define
$$\delta_3=min(1,\delta_1),\quad \delta_4=max(1,\delta_2).$$
From (\ref{NN}) and the boundedness given by (\ref{bcon}), it is observed that 
\begin{equation}
\label{bdd8}
\begin{aligned}
\delta_3\leq&\Vert D\bd{\mathcal{N}}(\bd{w}_2)\Vert\leq\delta_4 \\
\delta_3\Vert\bd{w}_1\Vert_{\mathcal{X}}^2\leq&(\bd{w}_1,D\bd{\mathcal{N}}(\bd{w}_2)\bd{w}_1)_{\mathcal{X}}\leq \delta_4\Vert\bd{w}_1\Vert_{\mathcal{X}}^2
\end{aligned}
\end{equation}
for every $\bd{w}_1,\bd{w}_2\in\mathcal{X}$, and thus Assumption \ref{gasp1} is satisfied. 
Furthermore, from definition of $\bd{\mathcal{A}}$ and $\bd{\mathcal{N}}(\cdot)$,
\begin{equation}
\label{bdd4p}
\begin{aligned}
(\bd{w},\bd{\mathcal{A}}\bd{\mathcal{N}}(\bd{w}))_{\mathcal{X}}\geq \lambda\Vert\bd{w}\Vert^2_{\mathcal{X}}
\end{aligned}
\end{equation}
for $\bd{w}\in\mathcal{D}(\bd{\mathcal{A}}^{1/2})$ such that $\bd{\mathcal{N}}(\bd{w})\in\mathcal{D}(\bd{\mathcal{A}})$; thus, Assumption \ref{gasp2} is satisfied.

Finally, the nonlinear operator $\bd{\mathcal{R}}(\cdot,t)$ is a composition of smooth functions of the potential vector $\bd{\varphi}$ and the vector $[\text{sat}_y(c_2),\ldots,\text{sat}_y(c_4)]^T$. Furthermore, $\bd{\varphi}$ is a Fr\'echet differentible function of $[c_1,c_2,\ldots,c_4]^T$. It is also observed that the variation of $\bd{\varphi}$ and $\bd{c}$ are bounded by the implication of the saturation functions $\text{sat}(\cdot)$ and $\text{sat}_y(\cdot)$ in (\ref{bal2-orig}); thus, $\bd{\mathcal{R}}(\cdot,t)$ is Lipschitz continuous with respect to $\bd{c}$; in other words, the nonlinearity of the system satisfies Assumption \ref{asp3}. Finally, the input vector $\bd{u}(t)$ is assumed to satisfy Assumption \ref{asp1}. The proof is then completed by Theorem \ref{ABDD}.
\eop

From Corollary \ref{LM}, eigenfunctions of $\bd{\mathcal{A}}$ can be used to approximate the system such that a subsequence of the approximate solutions converges to a solution of (\ref{bal2-orig}).
For the sake of simplicity, since the electrolyte concentration does not experience much change along the cell in time, it is set to be constant as in \cite{mp-r16} to find the eigenfunctions.  
For the solid concentration, $c_2$-$c_4$, the eigenfunctions are derived from the following eigenvalue problems: for $k=1,\ldots,3$,
\begin{align}
\label{EV}
\frac{1}{r_k^2}\frac{\partial }{\partial r_k}( r_k^2\frac{\partial \bd{z}_{k+1}}{\partial r_k} )=\lambda_{k}\bd{v}_{k+1}
\end{align}
in which the linear operator's domain is defined in (\ref{dom}). Solving (\ref{EV}) leads to finding the eigenfunctions as
\begin{equation}
\label{EVEF}
\begin{aligned}
  \bd{v}_{k+1}=
  \left\{ 
  \begin{array}{l l}
    1 & \quad 			\text{if $j=0$}\\
    \frac{\sin\big(\frac{\gamma_j}{R_k}r_k\big)}{r_k} & \quad \text{otherwise}\\
  \end{array} \right.
\end{aligned}
\end{equation} 
where $\gamma_j$ satisfies
$$\gamma_j=\tan(\gamma_j).$$

Note that, in the original electrochemical equations the derivatives of the solid concentrations $(c_2,c_3,c_4)$ with respect to the spatial variable $x$ are not involved. In order to add more accuracy to the system's solution, in the next step, the electrolyte concentration, $c_1$, is approximated by a piece-wise linear function instead of a constant and included in the system's dynamics.

%The constraint equation (\ref{E212b}) is approximated by a finite-dimensional one. 
Linear spline functions are appropriate choices for approximating the potential vector $\bd{\varphi}$ since (\ref{E212b}) includes second order differentiation. The Galerkin method is then used to find finite-dimensional nonlinear approximate algebraic equations. 
%These equations are transformed into a dynamical form by applying time differentiation to both sides of the algebraic equations. 
%In this way, the finite-dimensional approximation of the operator (\ref{IMP}) is automatically achieved as a part of the finite-dimensional dynamical form. 
%======================================================================
\section{Simulations and comparison to experimental data}
\begin{figure}

\centering
\subfloat[Charging ]{
  \includegraphics[width=60mm]{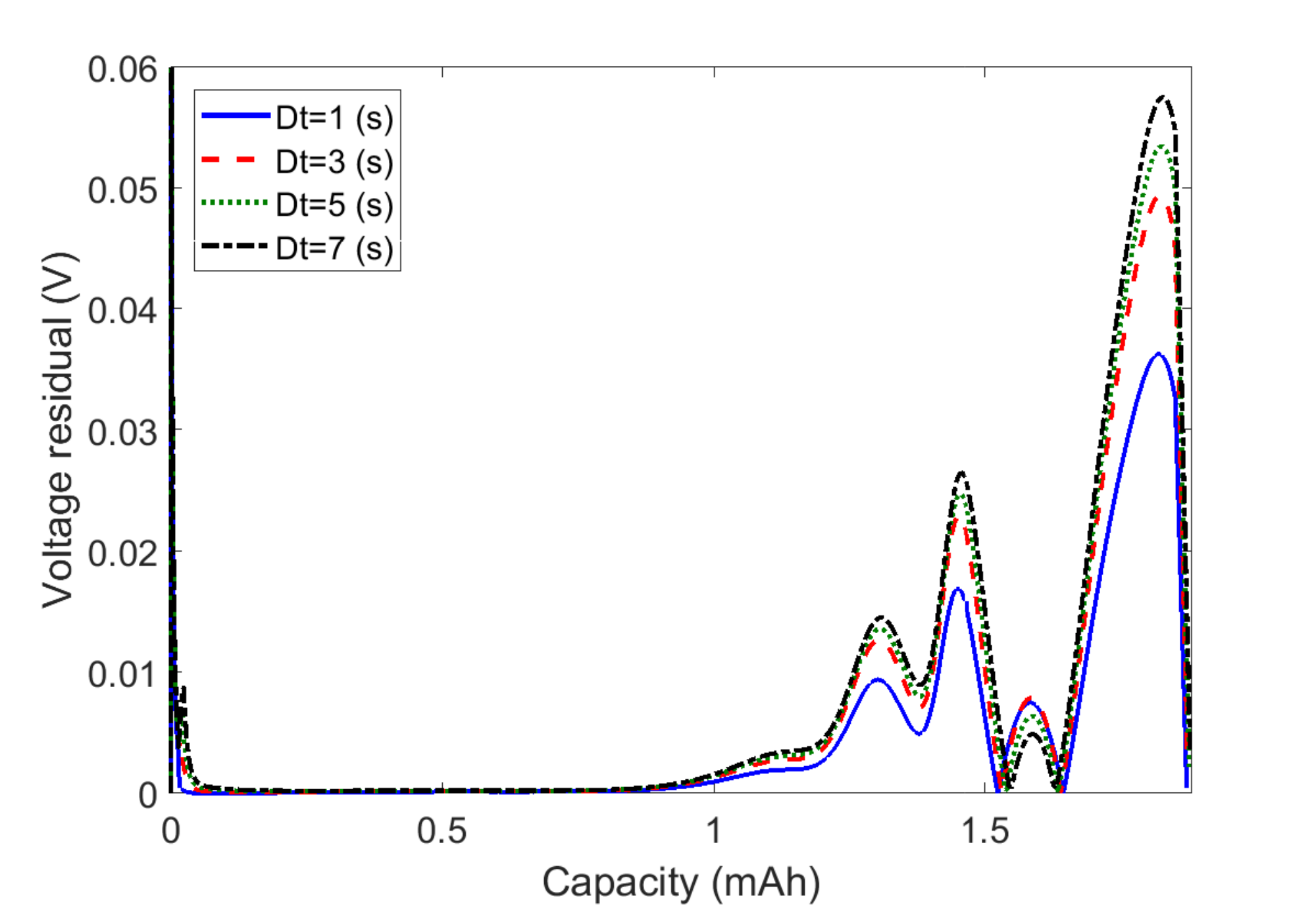}
}
\subfloat[Discharging ]{
  \includegraphics[width=60mm]{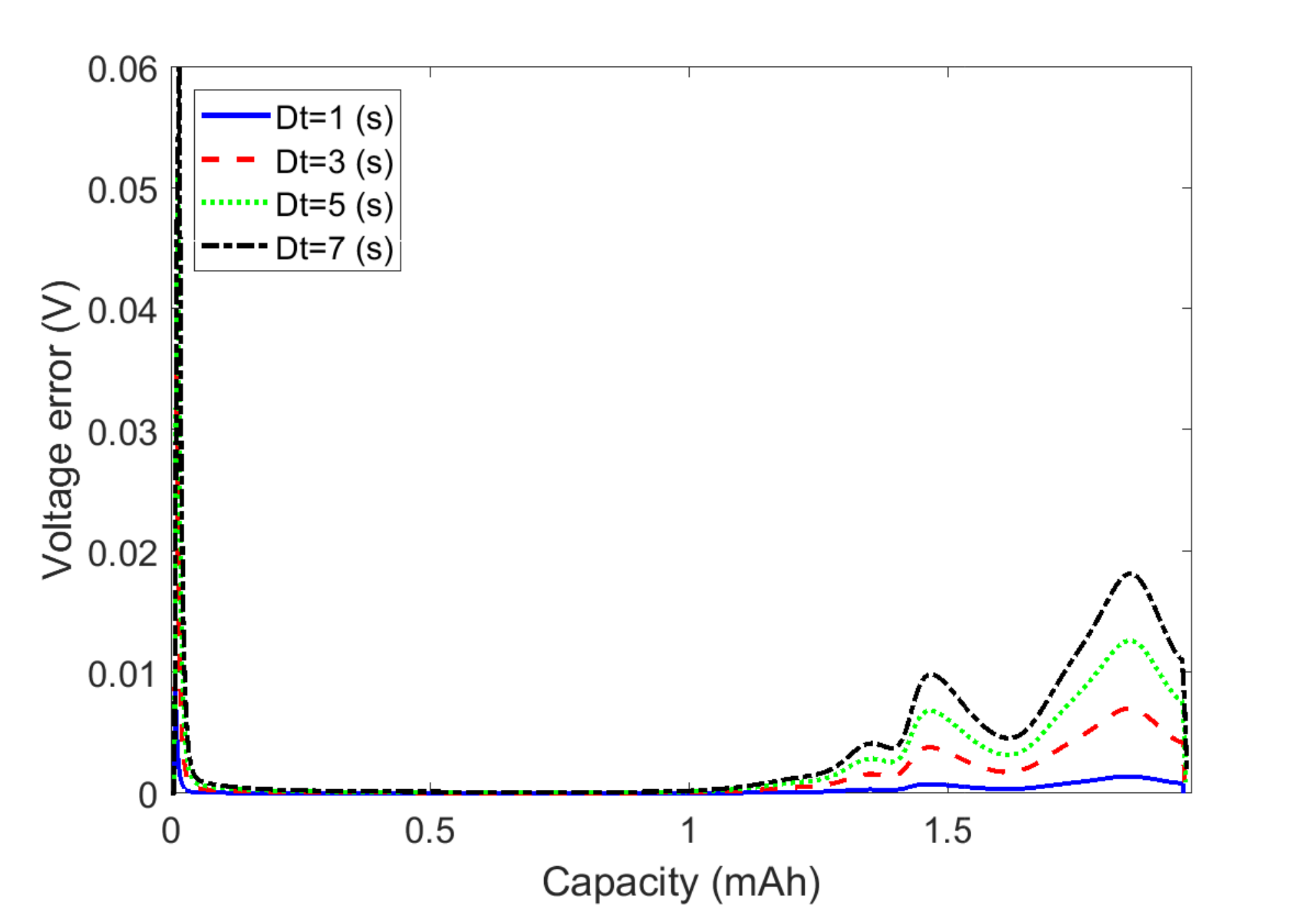}
}
\caption{Residual voltage for different correction sampling time $Dt$ compared to a reference voltage with the correction sampling time $Dt=0.5\,\text{s}$ at the charging/discharging current rate $1\text{C}$; the simulations used $N_1=4$ in the separator domain, $N_2=4$ along the positive electrode, and $N_3=6$ for every particle. It is observable that $Dt=3\, \text{s}$ provides a small residual voltage.}
\label{F2EVCR}
\end{figure}

\begin{figure}
\centering
\subfloat[Charging ]{
  \includegraphics[width=60mm]{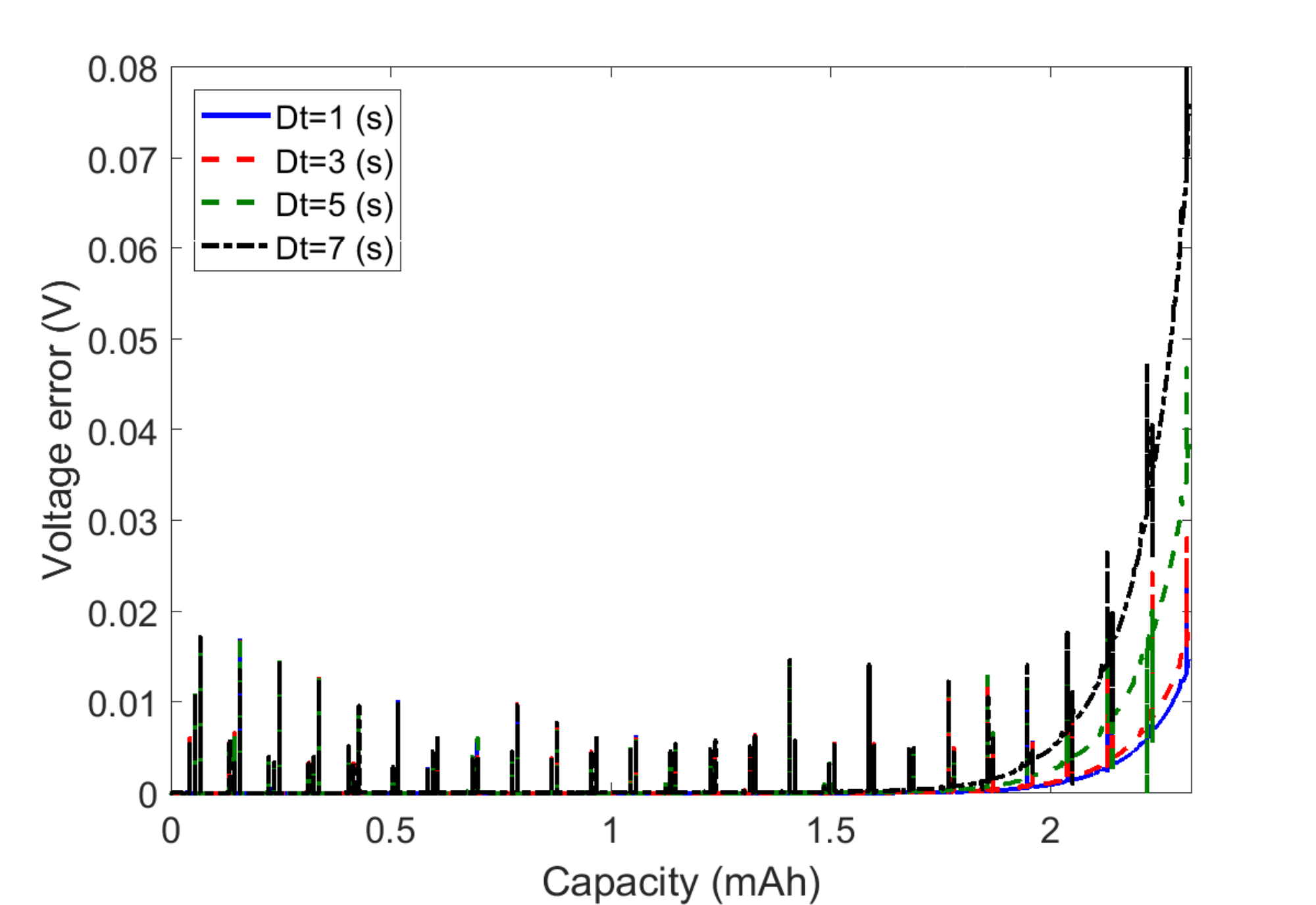}
}
\subfloat[Discharging ]{
  \includegraphics[width=60mm]{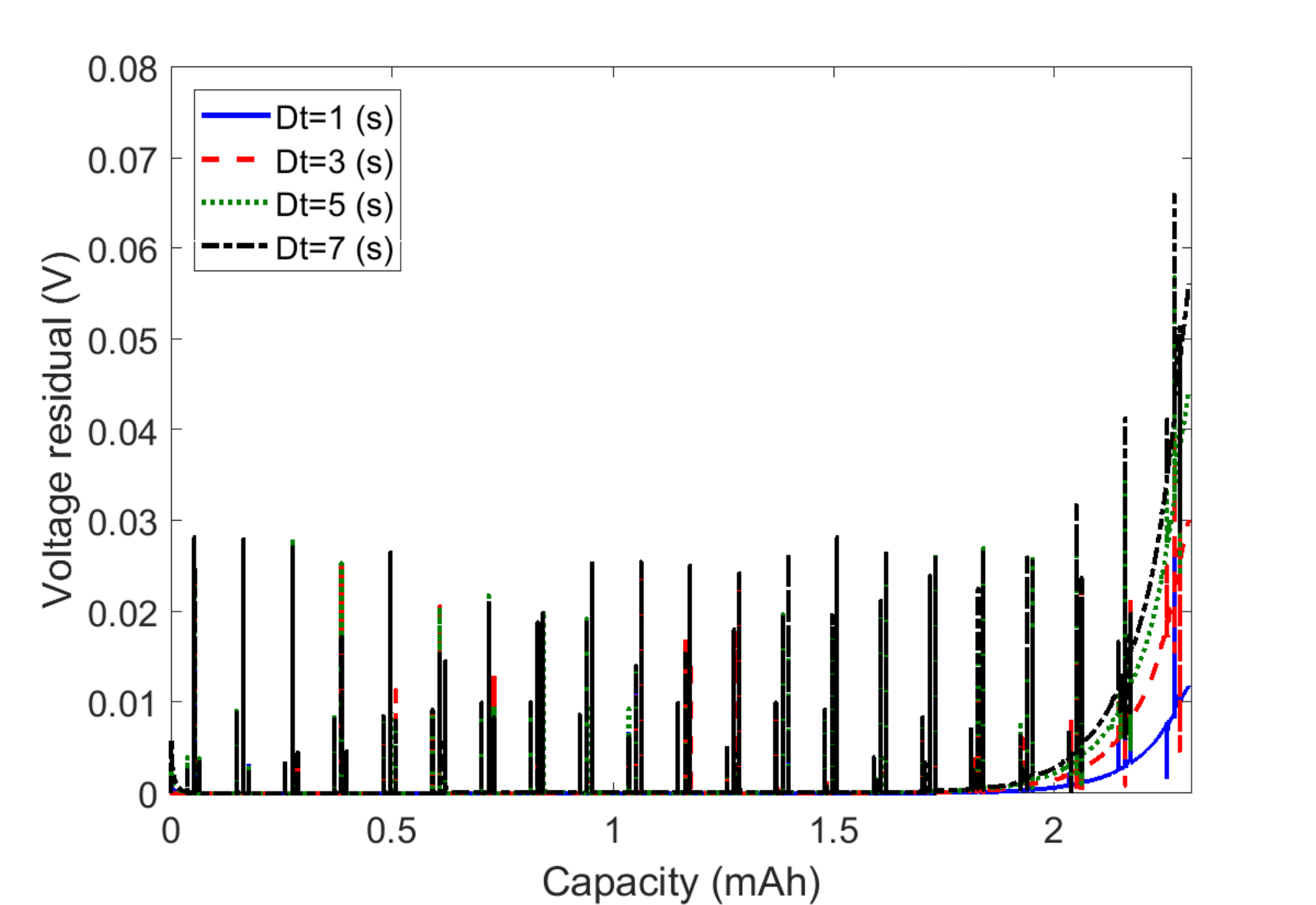}
}
\caption{Residual voltage for different correction sampling time $Dt$ compared to a reference voltage with the correction sampling time $Dt=0.5\, \text{s}$ at the charging/discharging impulsive current; the simulations used $N_1=4$ in the separator domain, $N_2=4$ along the positive electrode, and $N_3=6$ for every particle. It is observable that $Dt=3\, \text{s}$ provides a small residual voltage.}
\label{F2EVIM}
\end{figure}

The finite-dimensional approximation of the original electrochemical equations is a system differential algebraic equations (DAEs). However, using  the fully dynamical form (\ref{CDNA-orig}) leads to a system of ordinary differential equations (ODEs). 
Using time differentiation to convert DAEs into ODEs can introduce  inaccuracy in the form of an accumulation error. The accuracy of this approach was improved by periodically solving the constraint equations for the potential vector. The solutions for different sample periods $Dt$ are compared in Figure \ref{F2EVCR} and Figure \ref{F2EVIM} by showing their difference from a solution obtained for $Dt=0.5\,\text{s}$. It is observed that the difference between solutions approaches zero.

\begin{figure}
\centering
\subfloat[Charging current rate=0.2C]{
  \includegraphics[width=65mm]{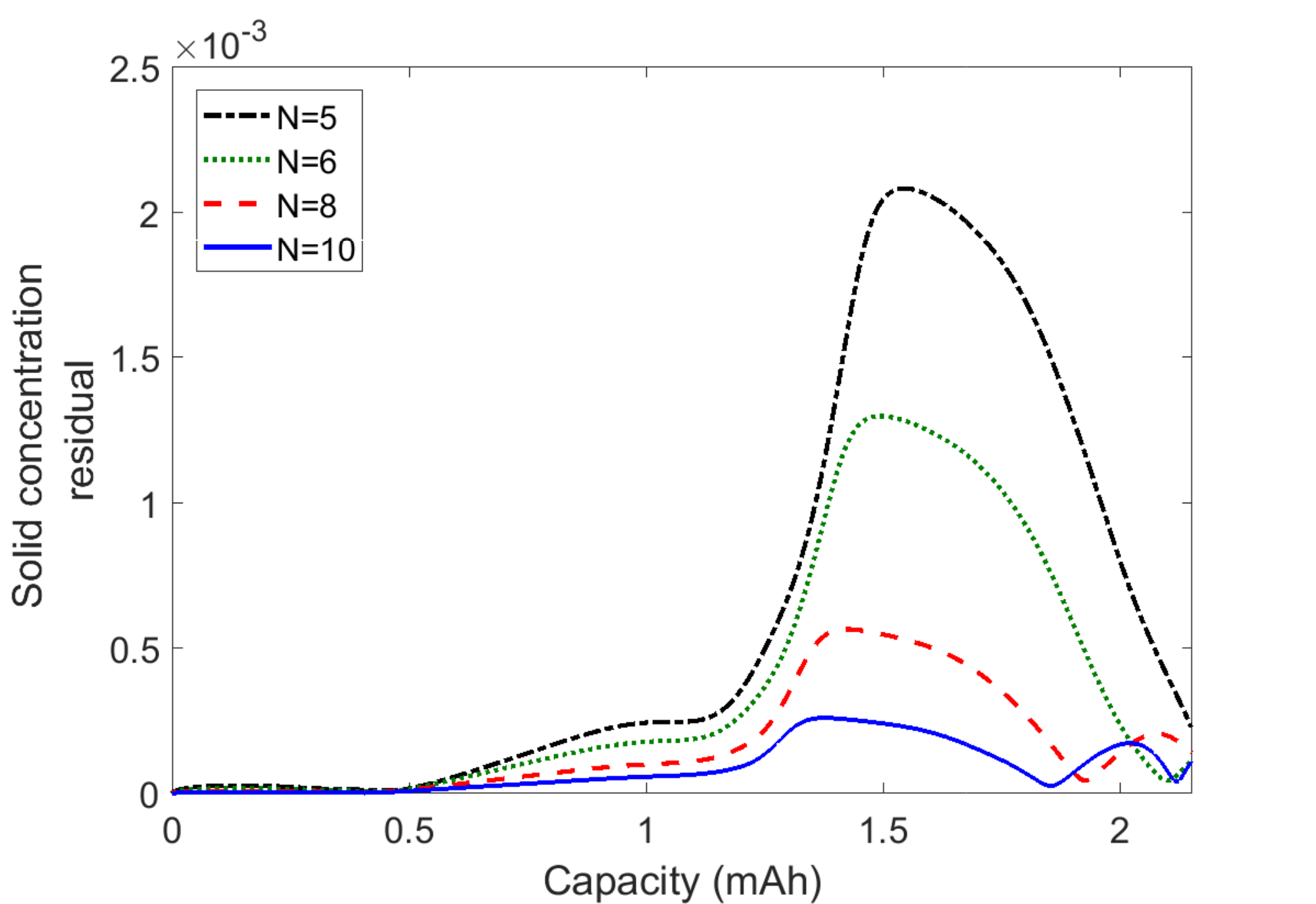}
}
\subfloat[Charging current rate=0.5C]{
  \includegraphics[width=65mm]{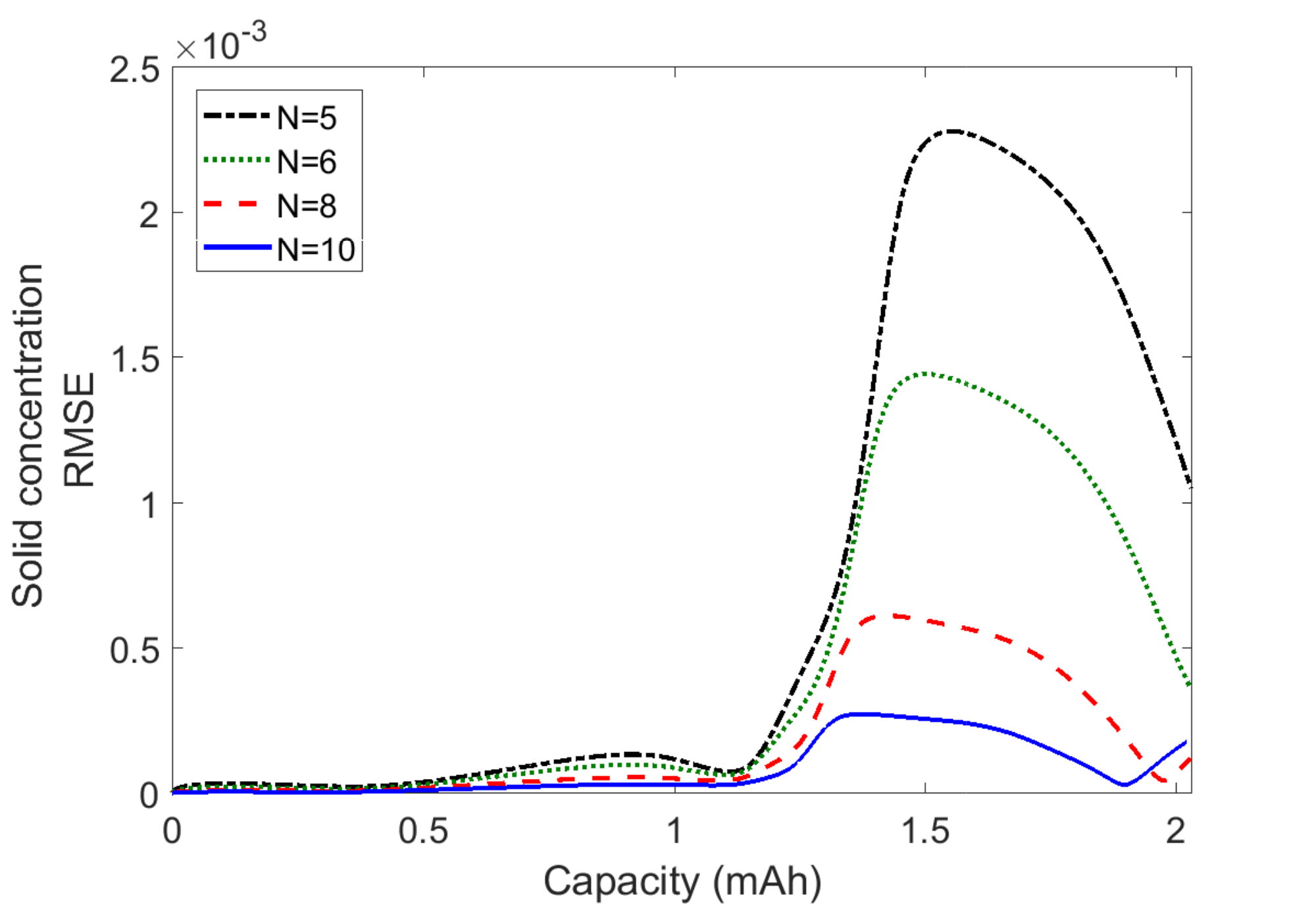}
}
\hspace{0mm}
\subfloat[Charging current rate=1C]{
  \includegraphics[width=65mm]{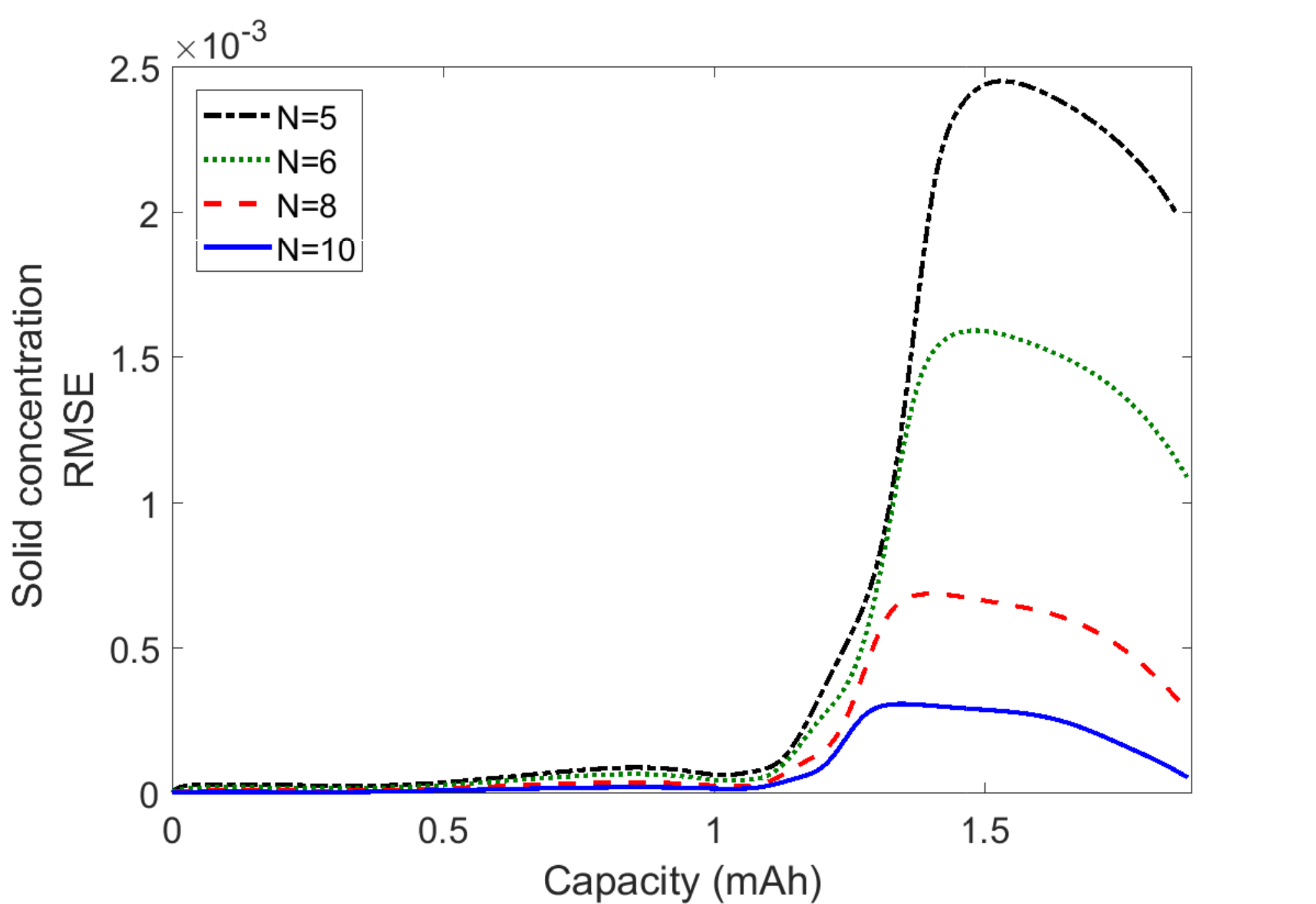}
}
\subfloat[Discharging current rate=0.2C]{
  \includegraphics[width=65mm]{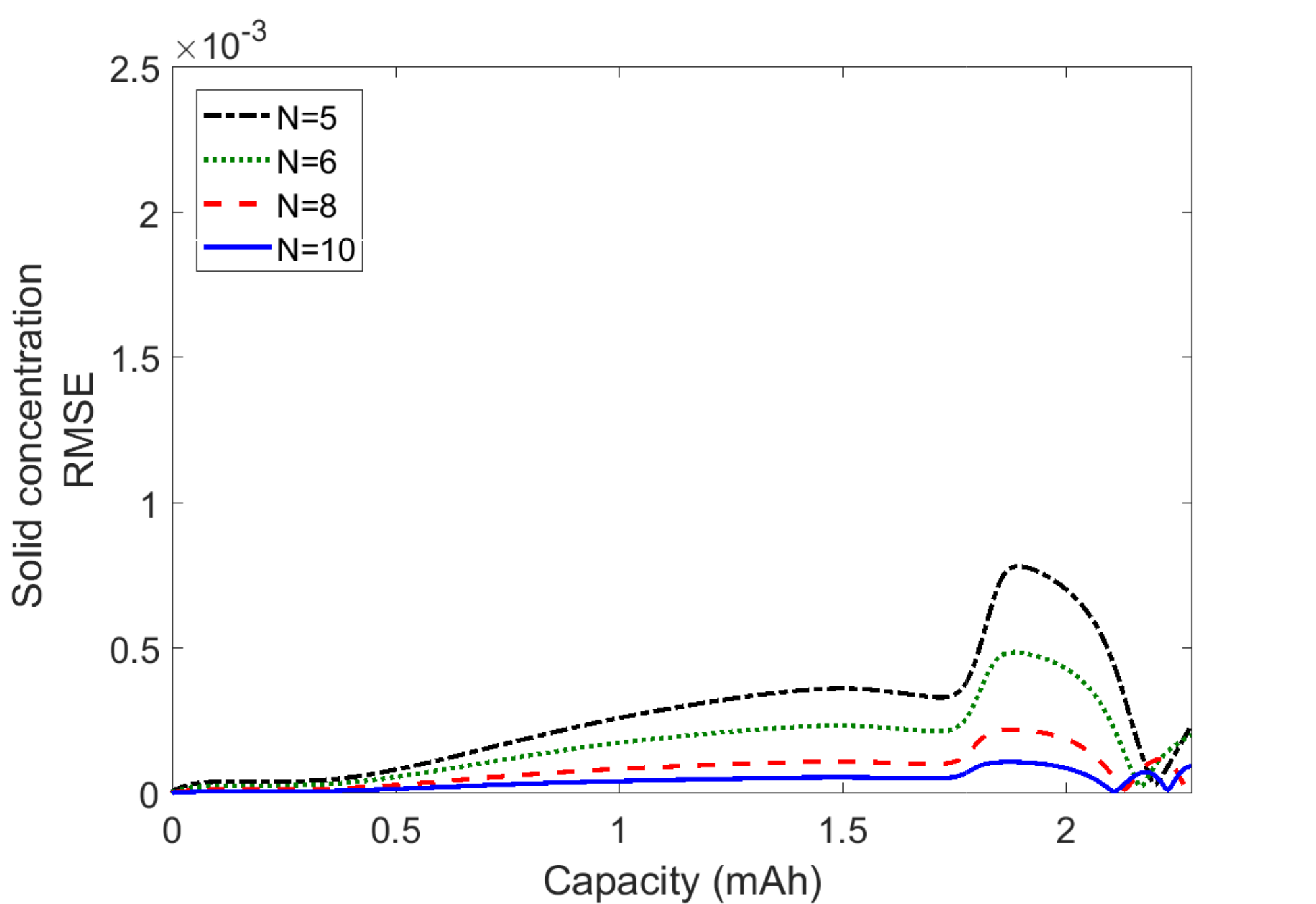}
}
\hspace{0mm}
\subfloat[Discharging current rate=0.5C]{   
  \includegraphics[width=65mm]{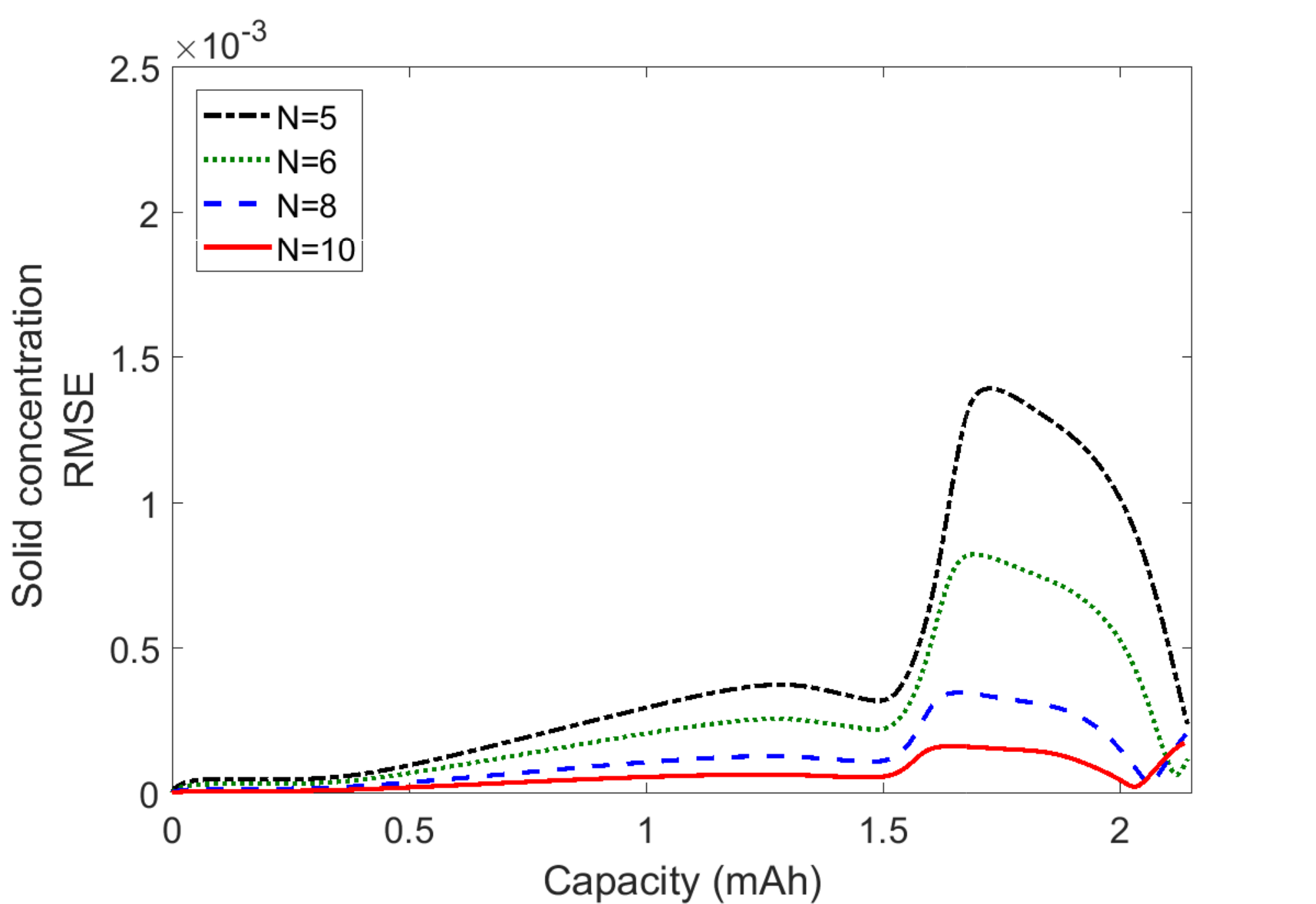}
}
\subfloat[Discharging current rate=1C]{
  \includegraphics[width=65mm]{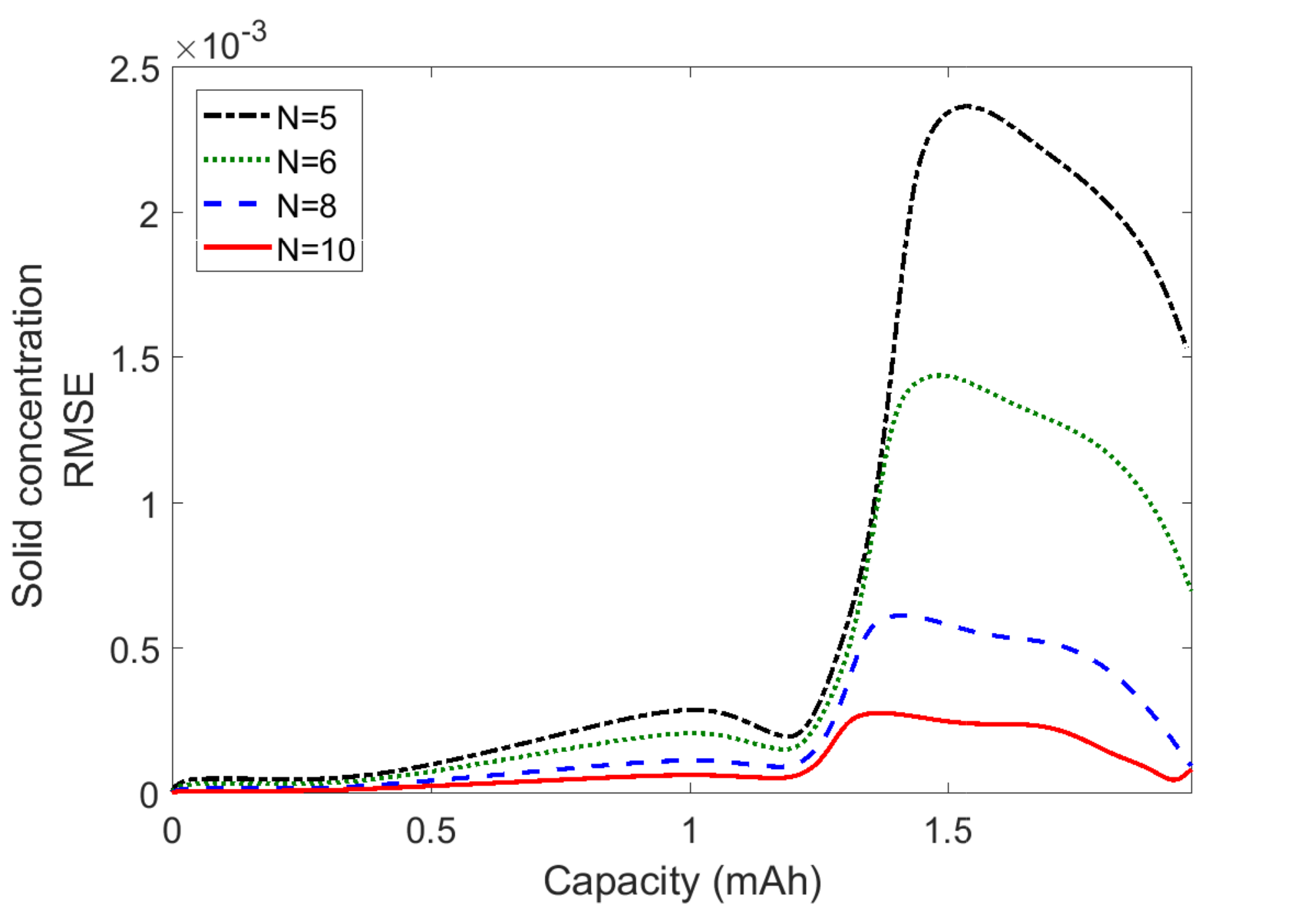}
}
\caption{RMSE between the solid concentration vector, at different orders of approximation, and the one with a large order of approximation $N_3=30$. The correction sampling time is $Dt=3\,\text{s}$. The approximate solutions converge fast especially at low current rates.}
\label{F7YE}
\end{figure}

\begin{figure}
\centering
\subfloat[Charging current rate=0.2C]{
  \includegraphics[width=65mm]{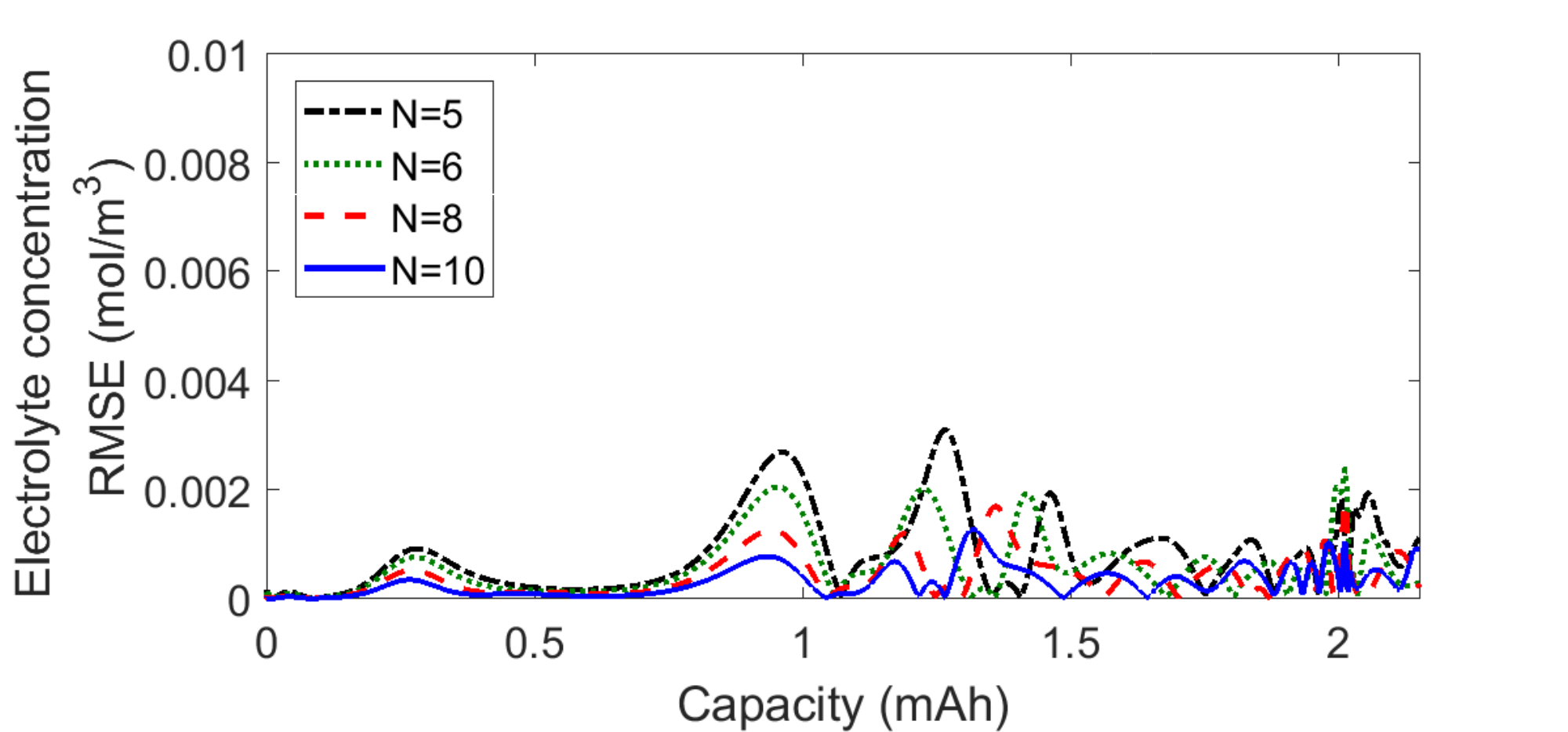}
}
\subfloat[Charging current rate=0.5C]{
  \includegraphics[width=65mm]{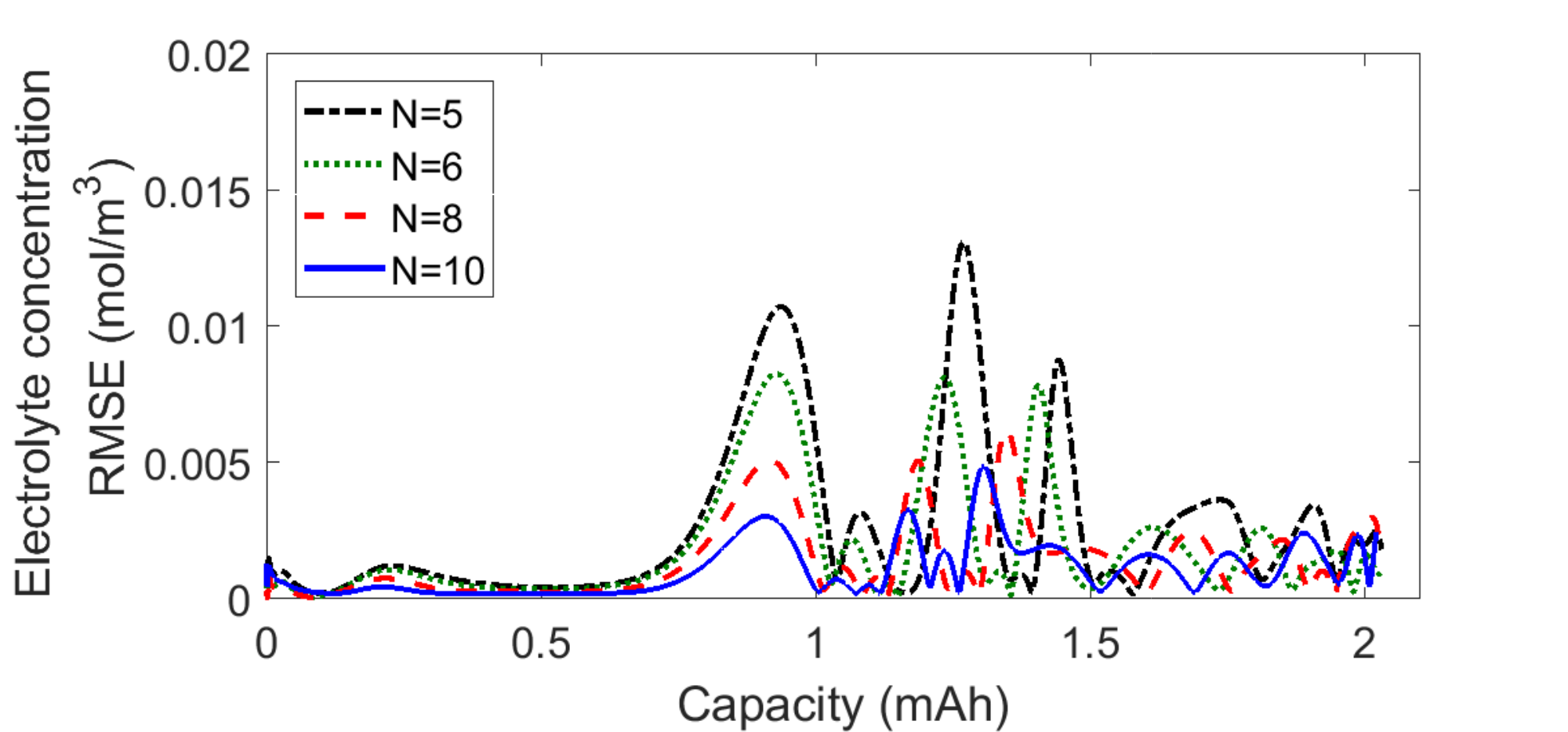}
}
\hspace{0mm}
\subfloat[Charging current rate=1C]{
  \includegraphics[width=65mm]{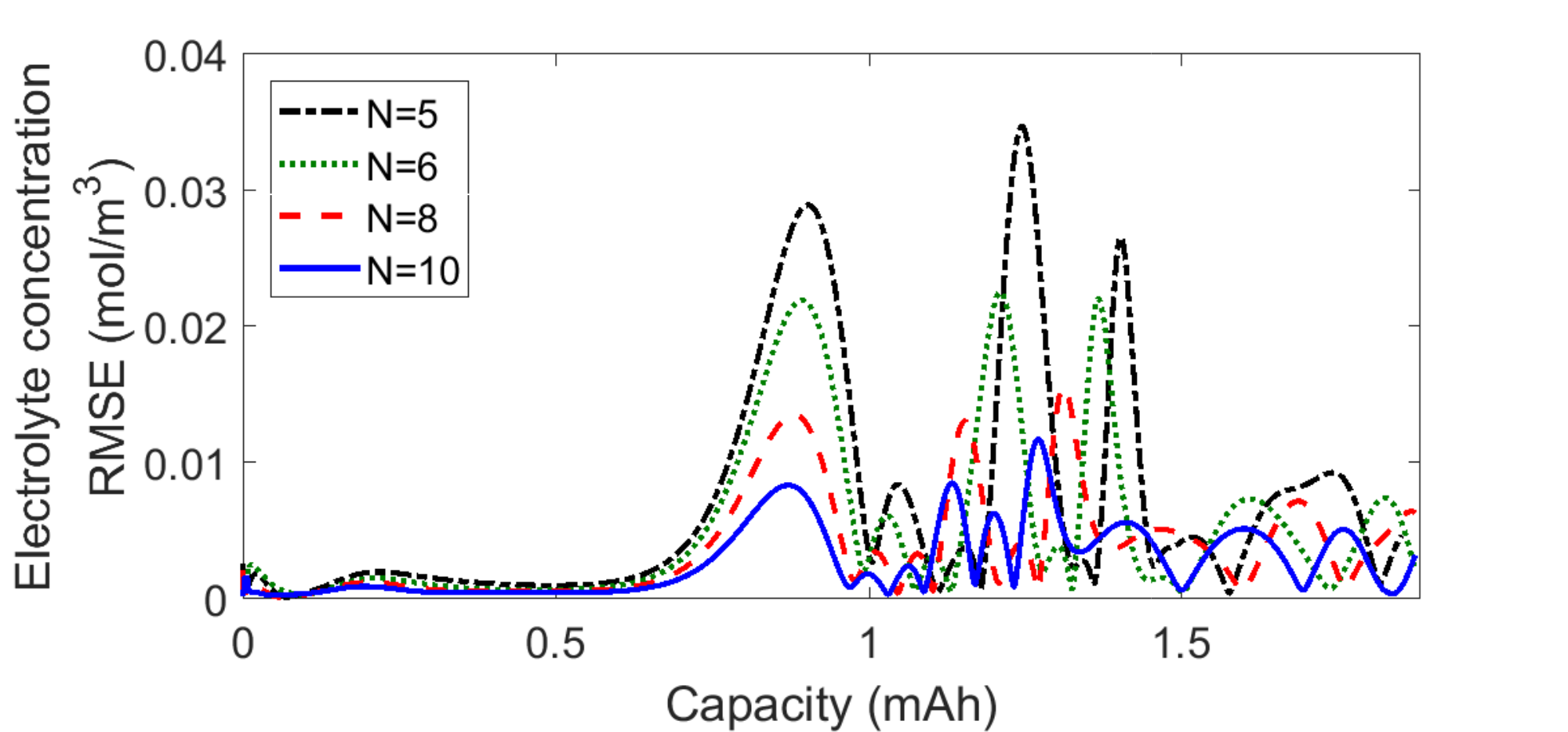}
}
\subfloat[Discharging current rate=0.2C]{
  \includegraphics[width=65mm]{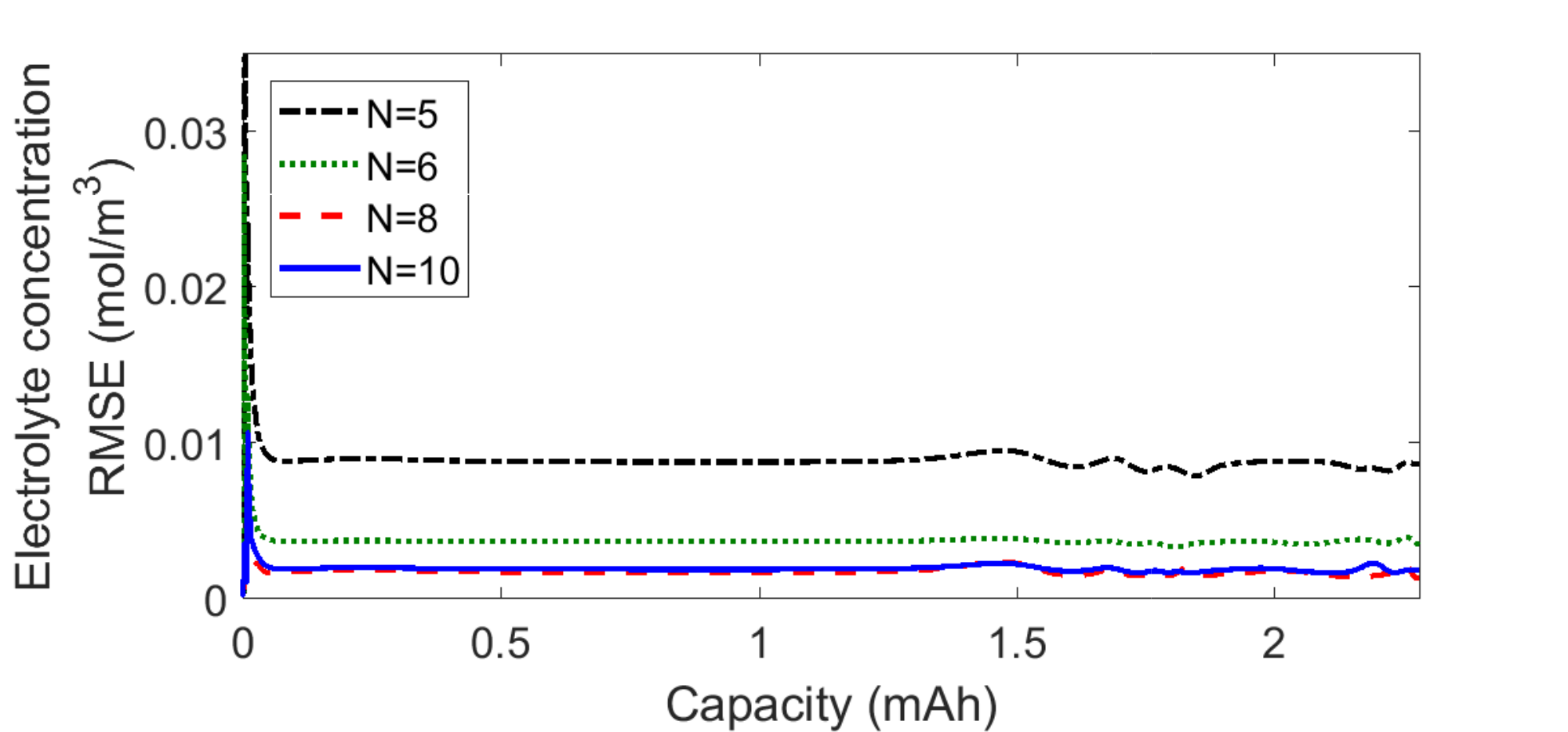}
}
\hspace{0mm}
\subfloat[Discharging current rate=0.5C]{   
  \includegraphics[width=65mm]{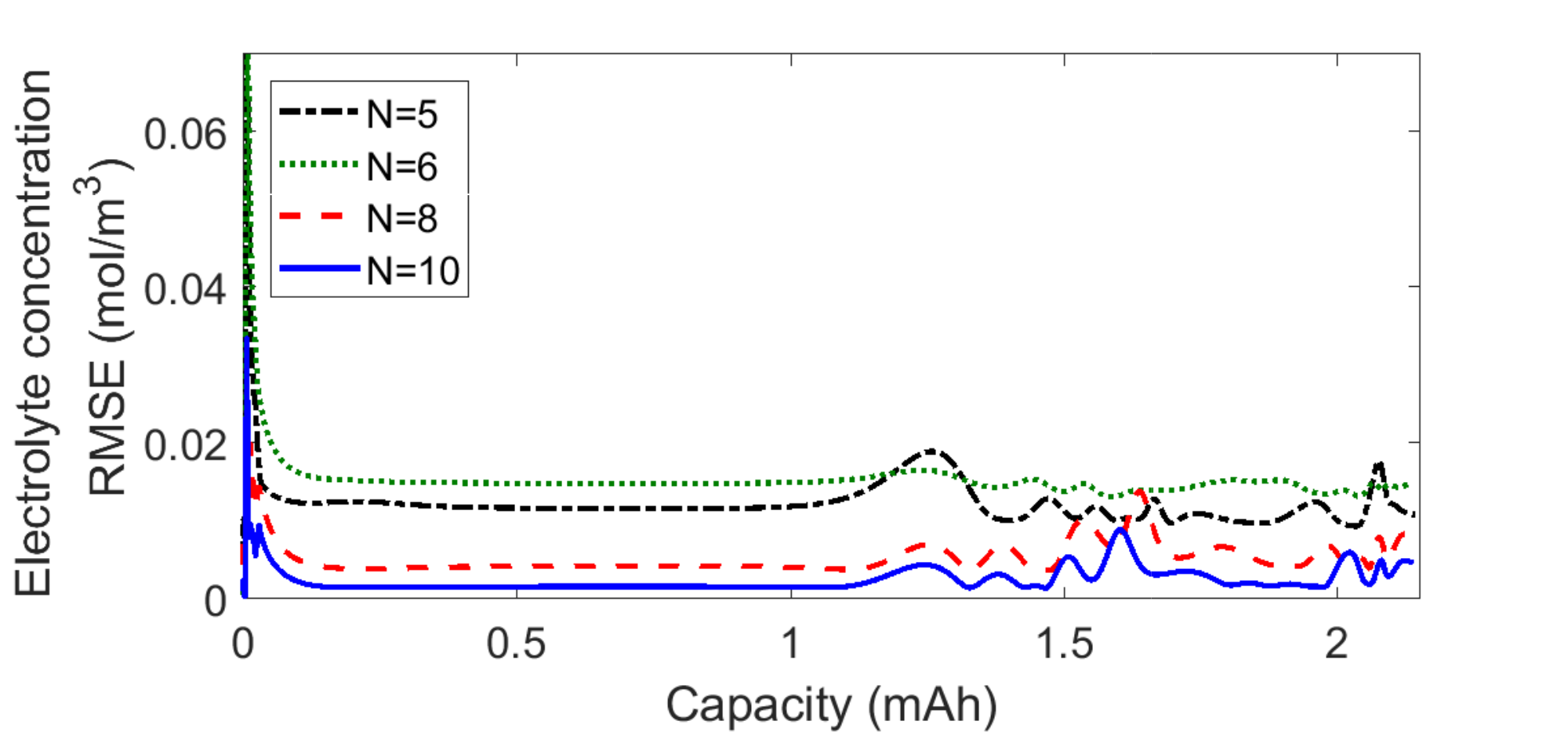}
}
\subfloat[Discharging current rate=1C]{
  \includegraphics[width=65mm]{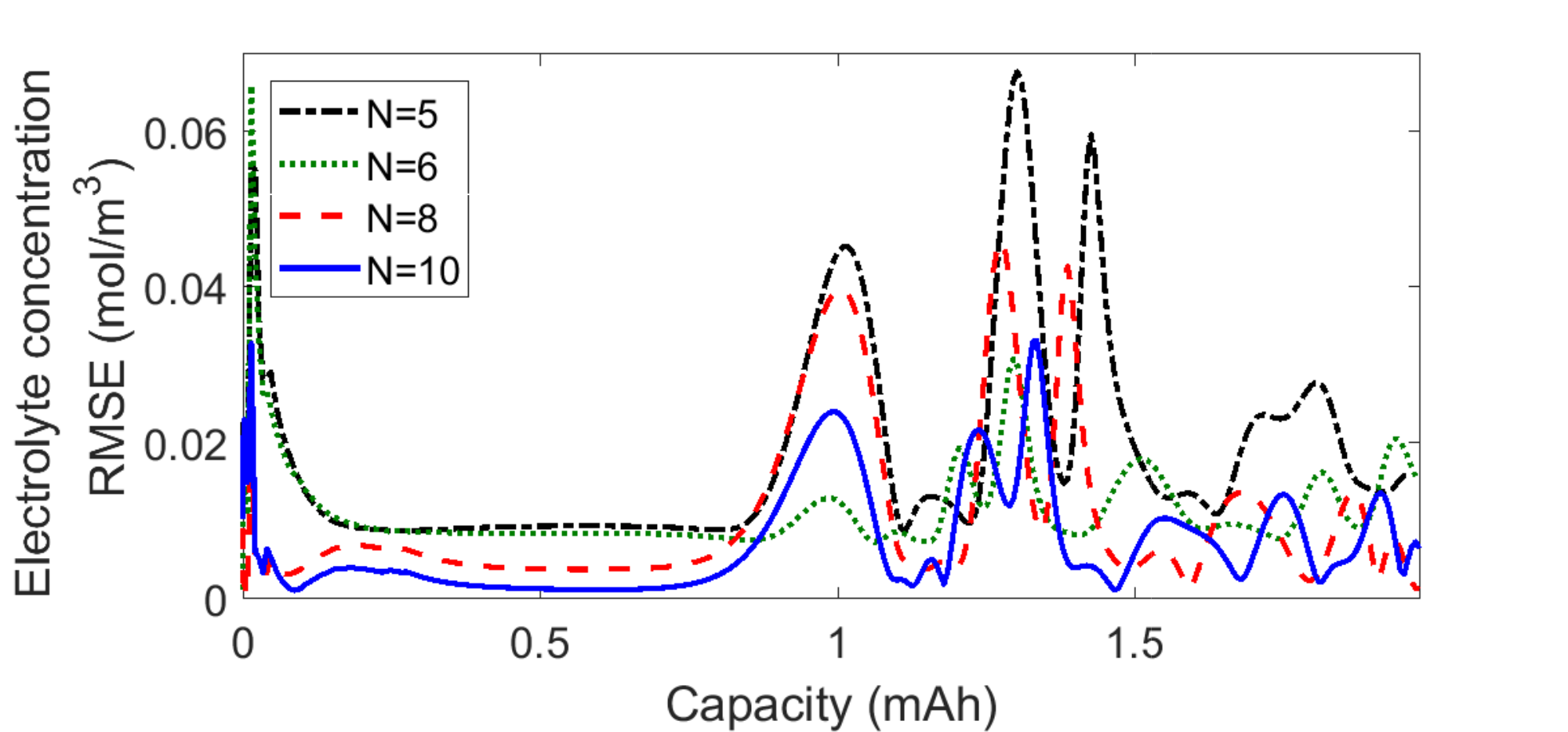}
}
\caption{RMSE between the electrolyte concentration vector, at different orders of approximation, and the one with a large order of approximation $N_3=30$. The correction sampling time is $Dt=3\,\text{s}$. The approximate solutions converge fast especially at low current rates.}
\label{F7CE}
\end{figure}

The convergence of solid and electrolyte concentrations as the order of approximation increases is shown by comparing the difference between the approximate solutions and a reference solution obtained by setting $N_3=30$, a large order of approximation. The root means square errors (RMSEs) between the solutions for solid and electrolyte concentration at different orders of approximation and the reference solution are shown in Figures \ref{F7YE} and \ref{F7CE}. The approximate solutions converge fast especially at low current rates.

The number of elements along the electrode and separator are denoted respectively by $N_1$ and $N_2$. The number of eigenfunctions along every particle is denoted by $N_3$. The simulations were run in MATLAB R2017  on a PC with Intel(R) CPU 2.3GHz processor and 32.0 GB RAM.

The experimental data used in this paper was generated in Laboratoire De R\'eactivit\'e Et Chimie Des Solides (LRCS) in Amiens, France.  In this experiment, the LFP electrode was recovered from a commercial graphite/LFP cell, LiFeBatt X2E (2.31 $\text{mAh}$, 40166, cell A) which is employed for hybrid electric applications \cite{mp-r16}. The cell underwent discharge to 2 $\text{V}$ at C/10 followed by a decrease of the current below C/50 while the potential was held to 2 $\text{V}$. Next, it was disassembled. Finally, the electrode whose area is 1.202 $\text{cm}^2$ was punched with a lithium metal foil for the counter electrode and a Whattman GF/D borosilicate glass fiber sheet for the separator to assemble a coin cell. For more details, please refer to \cite{mp-r16}.

\begin{figure} [ht]
\centering
\includegraphics[width=.6\linewidth]{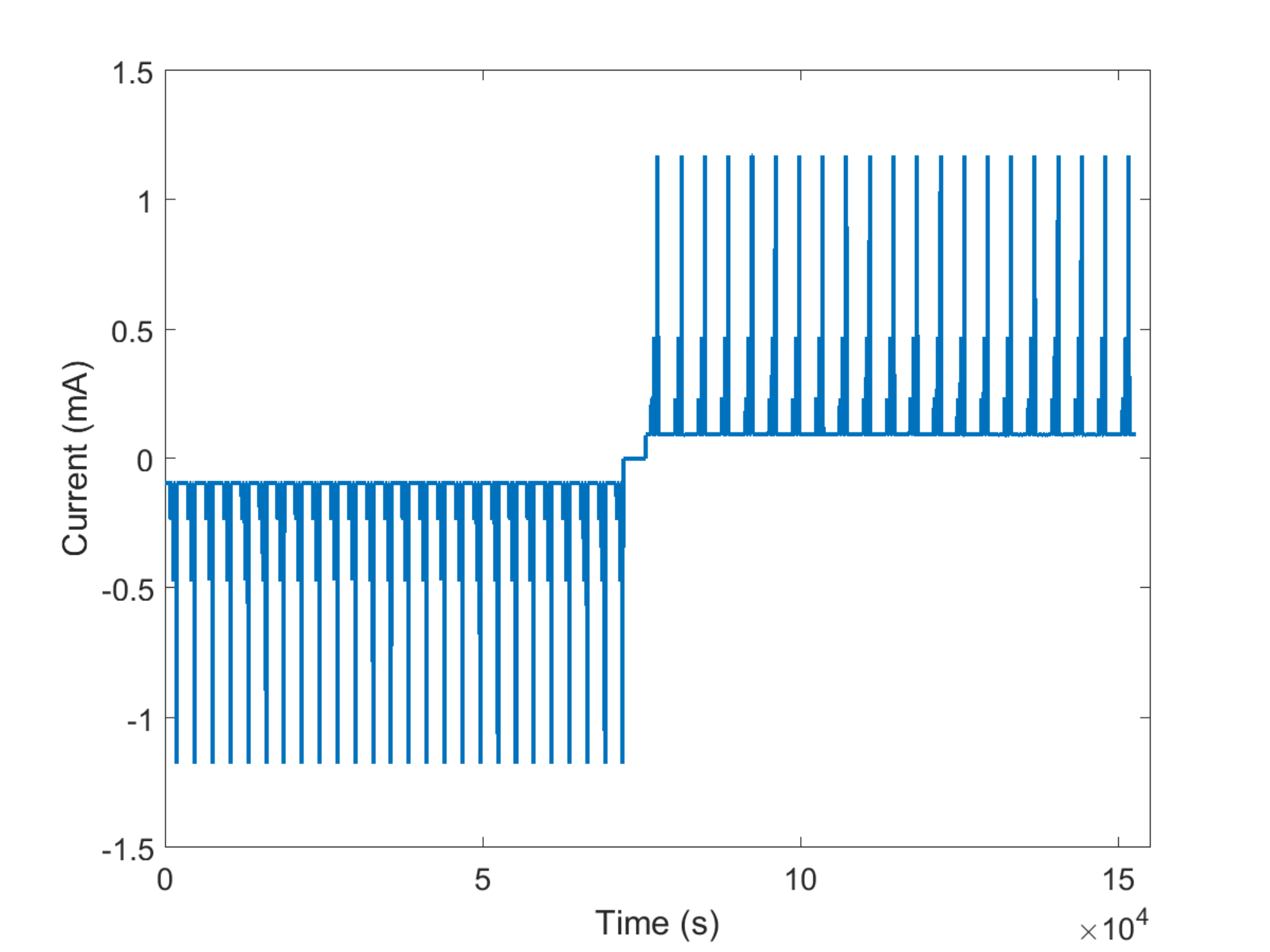}
\caption{Current profile in a charging and discharging cycle.}
\label{F1}
\end{figure}

\begin{figure}
\centering
\subfloat[Charging ]{
  \includegraphics[width=60mm]{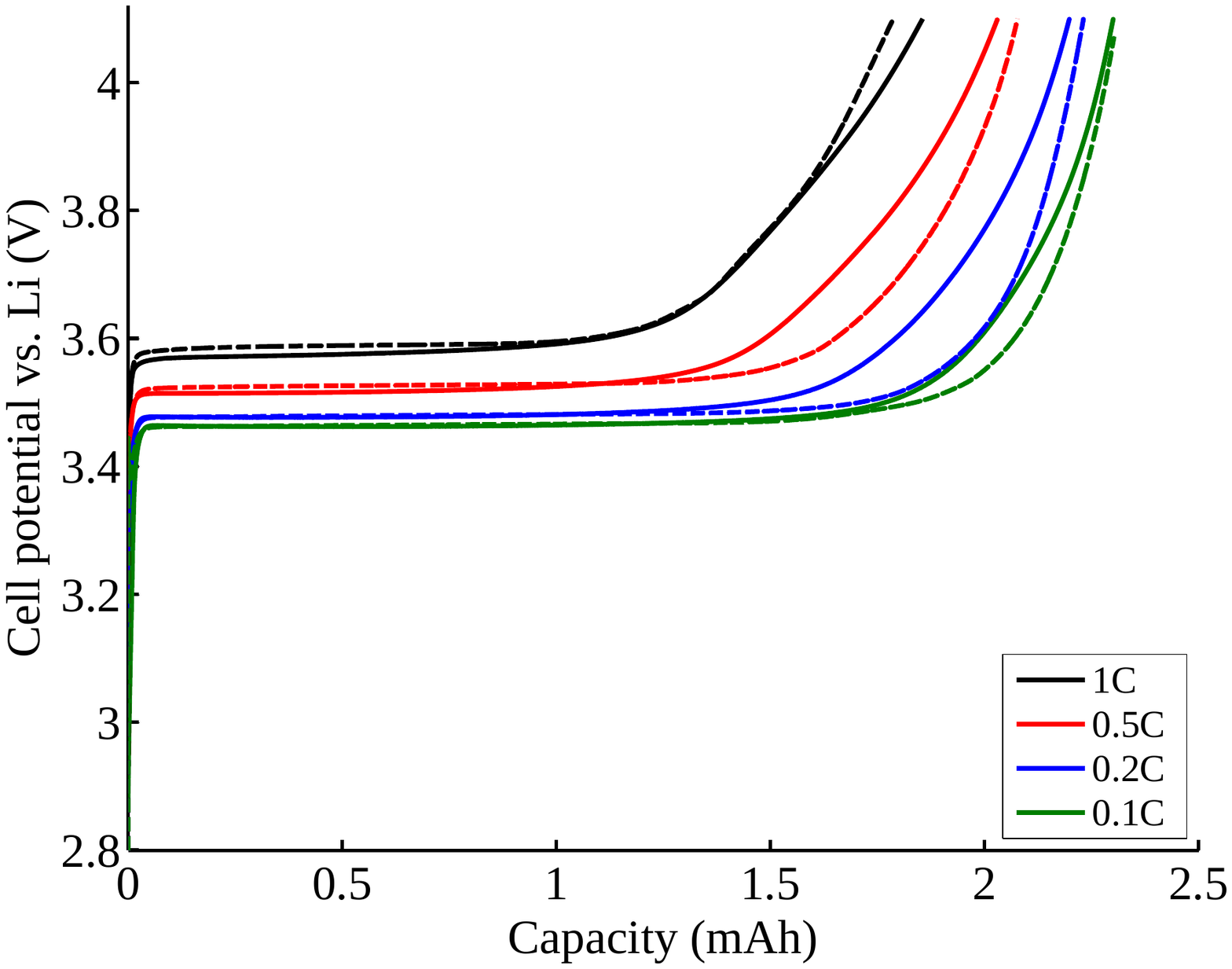}
}
\subfloat[Discharging ]{
  \includegraphics[width=60mm]{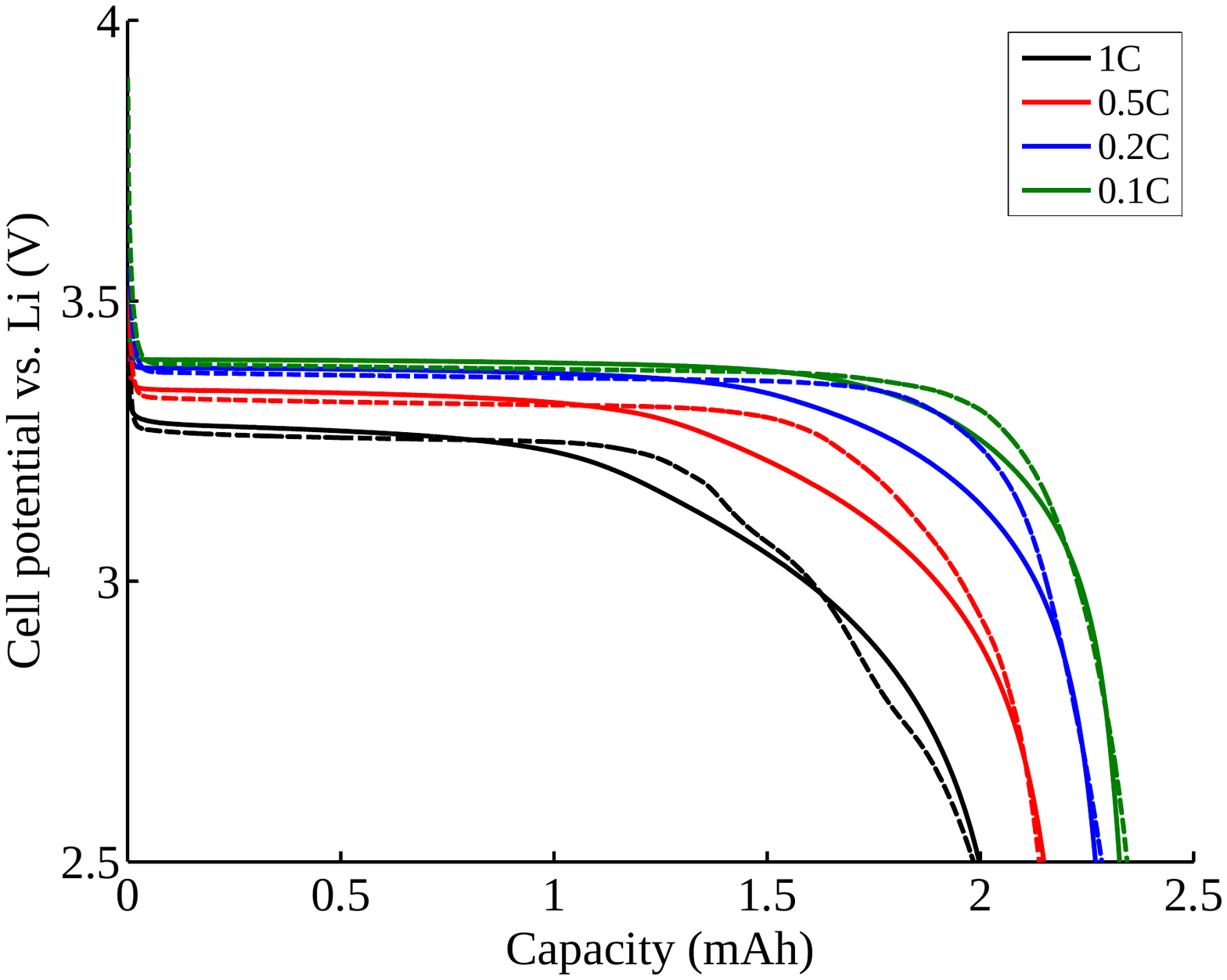}
}
\caption{Comparison of the simulation results of the fully dynamical representation with the experimental data for different current rates; the simulations used $N_1=4$ in the separator domain, $N_2=4$ along the positive electrode, $N_3=6$ for every particle, and the correction sampling time $Dt=3\, \text{s}$. A good agreement with the experimental data is observed. In these plots, the dashed and solid line respectively represent the simulation result and experimental data.}
\label{F2}
\end{figure}

\begin{figure} [ht]
\centering
\includegraphics[width=.6\linewidth]{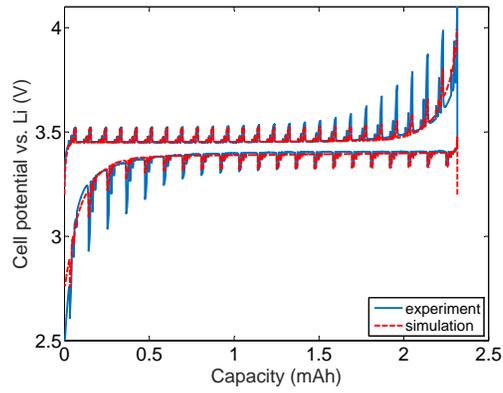}
\caption{Comparison of the simulation results of the fully dynamical representation with experimental data for the impulsive current; the simulations used $N_1=4$ in the separator domain, $N_2=4$ along the positive electrode, $N_3=6$ along every particle, and the correction sampling time $Dt=3\, \text{s}$. Agreement with the experimental data is observed. In these plots, the dashed and solid line respectively represents the simulation result and experimental data.}
\label{F4}
\end{figure}

Simulation results of solving the system equations  (\ref{bat5}), (\ref{bal2-orig})  for different charging and discharging current rates are shown in Figure \ref{F2} for $N_1=4$, $N_2=4$, and $N_3=6$ with correcting sampling time  $Dt=3\, \text{s}$. Similarly, the simulation results for an impulsive current (Figure \ref{F1}) are shown in Figure \ref{F4}. As observed from these figures  the results have a good agreement with  experimental data for constant current rates and most of the operation region with charging/discharging current profile presented in Figure \ref{F1}.

 The comparison of the computation time to the  experimental charging/discharging time is shown in Table \ref{T3}. The simulation time is   much faster than the actual time of the charging/discharging cycle.  The computation can be compared with the reported time in \cite{Coordinate} for the current rate $1\, \text{C}$. 
The reported MAPLE computation time in \cite{Coordinate} is using a 3.33 GHz Intel processor with 24 GB RAM for the degrees of freedom  $136$ and $72$ are respectively $28.361\,\text{s}$ and $9.812\,\text{s}$ which is comparable with $36\,\text{s}$ obtained in this paper for a degree of freedom $88$. The computational time is also less than the one introduced in \cite{urisanga2015}; the reported MAPLE computation time in \cite{urisanga2015} is $174.71\,\text{s}$. It should be noted that the nonlinearity of the equations in this paper is more than the nonlinearity  in either \cite{Coordinate} or \cite{urisanga2015}.

\begin{figure}
\centering
\subfloat[Charging ]{
  \includegraphics[width=60mm]{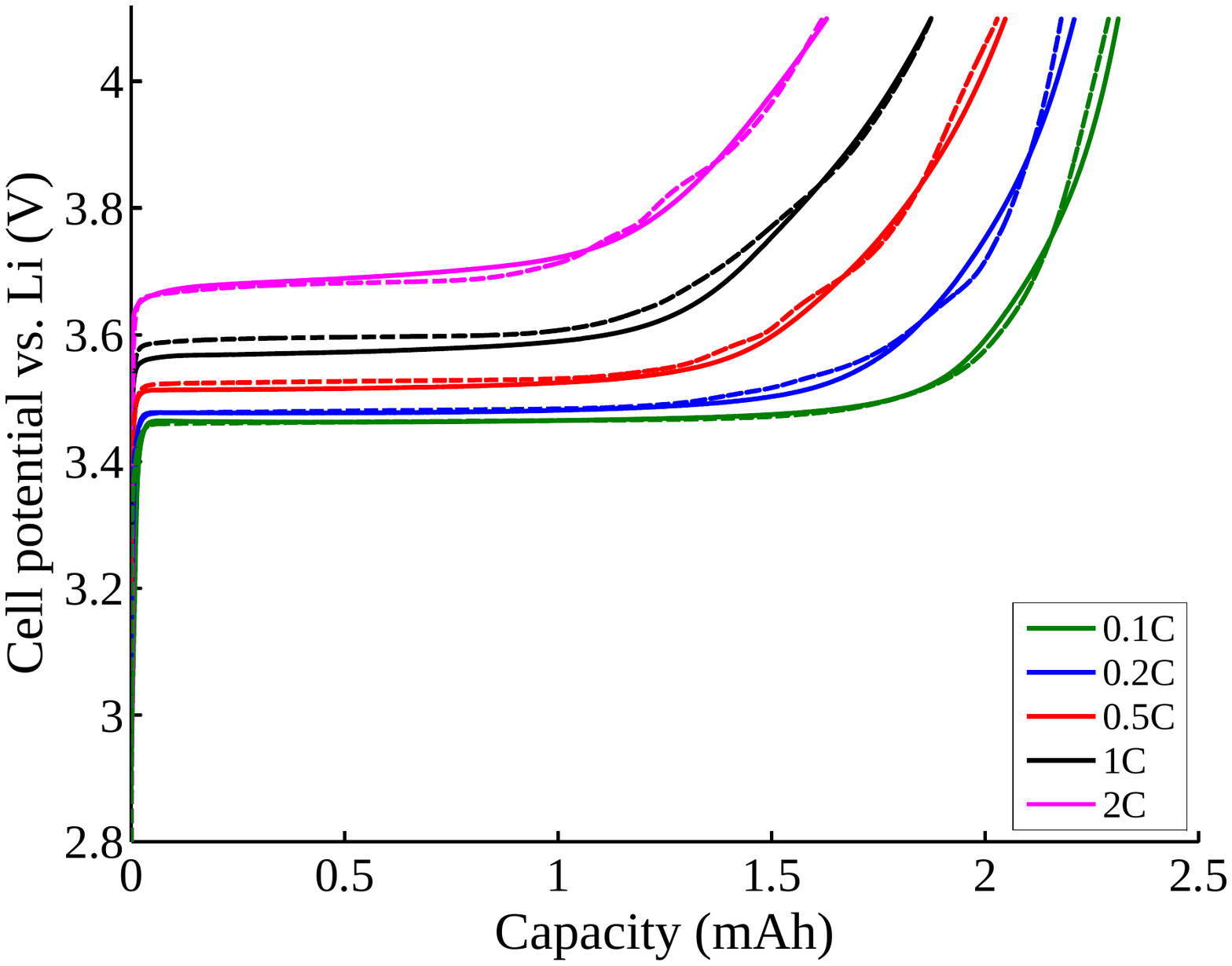}
}
\subfloat[Discharging ]{
  \includegraphics[width=60mm]{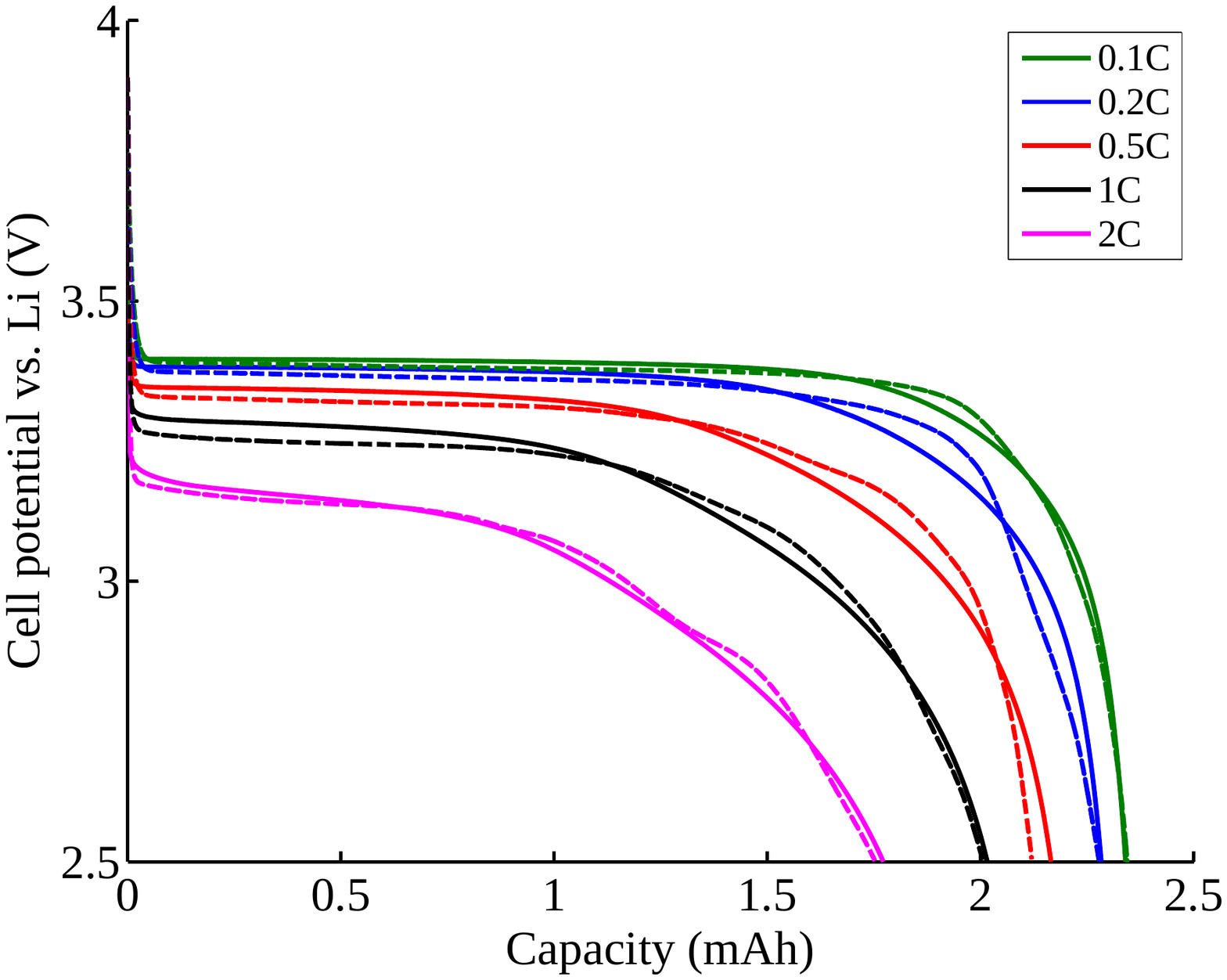}
}
\caption{Comparison of the simulations for equations with a rate dependent diffusion coefficient to the experimental data; the simulations used $N_1=4$ in the separator domain, $N_2=4$ along the positive electrode, $N_3=6$ along every particle, and correction sampling time $Dt=3\,\text{s}$. Agreement with the experimental data is improved compared to the standard variable solid-diffusivity model. In these plots, the dashed and solid line indicate the simulation result and experimental data respectively.}
\label{F9}
\end{figure}

A likely cause of the  discrepancy between the simulation results and experimental data  is errors in the modeling parameters. A more accurate model was obtained by including rate dependency in the diffusivity. This is done by changing the activity correction factor (\ref{acf}) to be rate dependent. For the discharging process,  the activity correction factor becomes
\begin{align*}
\alpha_k(y_k)&=9\exp(-25y_k)+15\omega_0 \exp(-30(1-y_k))\\
&+3\omega_1\exp(-15(1-y_k))+0.2\omega_2/(1+(y_k-0.5)^2)
\end{align*}
where $\omega_0$, $\omega_1$, and $\omega_2$ are rate dependent correcting coefficients.
Similarly, for charging process,  the activity correction factor becomes
\begin{align*}
\alpha_k(y_k)&= 9\omega_3\exp(-25y_k)+15 \exp(-30(1-y_k))+0.2\omega_4/(1+(y_k-0.5)^2)
\end{align*}
where $\omega_3$ and $\omega_4$ are rate dependent correcting coefficients (See Table \ref{T4} for the values of these coefficients at different current rates). 
The correcting coefficients are calculated to minimize the difference between the model output and the experimental data. The simulations with the modified activity correcting factor are shown in Figure \ref{F9}. It is observed that the rate dependent diffusion coefficient improves the match to the experimental data. 
\section{Conclusions}
In this paper, a model for lithium-ion cells that is accurate but also appropriate for real-time applications in hybrid vehicles was introduced. A challenging part of the real-time applications is obtaining a model that is precise and yet fast to be implemented online. As discussed before, in the case of LFP cells, the electrochemical models are accurate; however, they are composed of both constraint equations and dynamical equations.  Many simplified models, such as equivalent circuits,  are accurate for specified parameter values but extensions are difficult. The proposed model is computationally simple, but by being physics-based, it can be easily adjusted to different working conditions. 

Next, the approximation of the electrochemical equations of an LFP cell with a well-posed state space representation was considered. Unlike many simplified models introduced in the literature, the state space representation preserves most of the cell's dynamics. 
It was shown that the constrained equations are well-posed, and the solid and electrolyte potential were restated as  functions of the state vector $\bd{c}$.  In the next step, a nonlinear low-order approximation was developed based on the modes of the linear part of the model, which is known to preserve the key dynamical behaviour. 

Simulation results showed a good agreement with experimental data even for low-order of approximations. It was also observed from the simulations that the approximate solutions converge as the order of approximation increases. Furthermore, the simulation time was much faster than the time elapsed for the experiment. The computation time is also comparable with the one reported in the literature for solving the electrochemical equations with constant diffusion coefficient.

The introduced reduced order model in this paper has the same accuracy as the ones introduced in literature including \cite{tanim2015aging} for the dynamics with a constant diffusion coefficient. 
The model was further improved by including a rate dependent solid diffusion coefficient; the accuracy of the model was increased as a result of the rate-dependent variable diffusion coefficient.  In order to add more accuracy to the system modeling,  the OCV term could  be modeled by a dynamical model. In this way, the effect of hysteresis can be included in the modeling. 

Observer design for SOC estimation using the  fully dynamical model described in this paper is the object of current research.

% was obtained \cite{ } and showed good results. 

\section*{Acknowledgement}
The experimental data used in this paper was generated in Laboratoire De R\'eactivit\'e Et Chimie Des Solides (LRCS) in Amiens, France and the authors thank them for the  permission to use this data to compare with simulations. The financial support of Automotive Partnership Canada (APC), Ontario Research Fund (ORF), and General Motors for this research is gratefully acknowledged.
\section*{References}
%\bibliographystyle{elsarticle-harv}
%\bibliography{ref}

\newpage
\begin{table*}[t]
\caption{Lithium-ion cell  parameters. }
\label{T1}
\centering 
\resizebox{\textwidth}{!}{ 
\begin{tabular}{ p{1cm} c c c}
\hline \hline
parameter & definition & value (separator)& value (LFP electrode) \\
\hline
$l_{cat}$ & thickness of the negative electrode $(\text{m})$ & &$72\times 10^{-6} $\\
\hline
$l_{sep}$ & thickness of the separator $(\text{m})$ & $675\times 10^{-6}$&\\
\hline
$R_k$ & radius of the spherical solid particles $(\text{m})$ & &$  1.44\times 10^{-7},2.70\times 10^{-7},5.42\times 10^{-7}$\\
\hline
$R$ & gas constant $(\text{J} / (\text{mol}. \text{K}))$ & &$8.3145$\\
\hline
$F$ & Faraday's constant $(\text{A}. \text{s}/ \text{mol})$ & &$96485$\\
\hline
$t_+^0$ & transference number & $0.363$ &$ 0.363$\\
\hline
$\epsilon_e$ & volume fraction of the electrolyte phase & $ 0.6$&$ 0.5$\\
\hline
$k^{eff}$ & effective conductivity in the electrolyte phase $(\text{s}/\text{m})$ & $ 0.6042$&$ 0.4596$\\
\hline
$k^{eff}_D$ & $\frac{k^{eff}2RT(1-t^0_+)}{F}$ && \\
\hline
$\sigma^{eff}$ & effective conductivity in the solid phase $(\text{s} /\text{m})$ & &$  6.75$ \\
\hline
$D^{eff}_e$ & effective diffusivity in the electrolyte phase $(\text{m}^2/\text{s})$ & $ 4.028\times 10^{-10}$&$ 3.677\times 10^{-10}$\\
\hline
$\mathcal{D}$ & Diffusion coefficient of the spherical particle $(\text{m}^2 /\text{s})$ & &$ 4.21\times 10^{-18}$\\
\hline
$c_{s,max}$ & maximum solid state concentration $(\text{mol}/\text{m}^3)$ & &$ 22.860\times 10^{3}$\\
\hline
$i_0$ & exchange current density $(\text{A} /\text{m}^2)$ & $ 0$&$ 3.25\times 10^{-2}$\\
\hline

\end{tabular}
}
\end{table*}

\begin{table}
\caption{Filtering and Saturation functions parameters}
\label{T2}
\centering 
\begin{tabular}{c|c|c|c|c|c}
\hline
$g$ & $h_0$ & $h_1$ & $a_0$ & $b_0$ & $\epsilon_0$\\
\hline
$1$&$2$&$3$&$1$&$2.0251$&$.0001$\\
\hline
\end{tabular}
\end{table}

\begin{table}[t]
\caption{Simulation and experiment charging/discharging process time. }
\label{T3}
\centering 
\begin{tabular}{ c c c c c c}
\hline \hline
Input current & $0.1\text{C}$ & $0.2\text{C}$ & $0.5\text{C}$ & $1\text{C}$ & Impulsive\\
\hline
\multicolumn{6}{c}{Computation time for charging cycle (s)}\\
\hline
experiment &$35859$&$17135$&$6330$&$2892$&$13447$\\
simulation &$479$&$251$&$76$&$33$&$361$\\
\hline
\multicolumn{6}{c}{Computation time for discharging cycle (s)}\\
\hline
experiment &$35785$&$17448$&$6620$&$3071$&$13407$\\
simulation &$453$&$212$&$98$&$36$&$356$\\
\hline

\end{tabular}\\
\end{table} 

\begin{table}
\caption{Diffusion coefficient correction factors.}
\label{T4}
\centering 
\begin{tabular}{c|c|c|c|c|c}
Current rate& $2\text{C}$&$1\text{C}$&$0.5\text{C}$&$0.2\text{C}$&$0.1\text{C}$\\
\hline
$\omega_0$ & $1.1$ & $0.9$ & $1$ & $1.2$ & $1.3$\\
\hline
$\omega_1$ & $1.2$ & $0.96$ & $0.6$ & $0.4$ & $0.3$\\
\hline
$\omega_2$ & $1$ & $0.8$ & $0.5$ & $0.3$ & $0.3$\\
\hline
$\omega_3$ & $1.5$ & $1$ & $0.75$ & $0.85$ & $1$\\
\hline
$\omega_4$ & $1.95$ & $1.3$ & $0.75$ & $0.35$ & $0.4$\\
\end{tabular}
\end{table}

\begin{figure}
\centering
\subfloat[Charging current rate=0.2C]{
  \includegraphics[width=65mm]{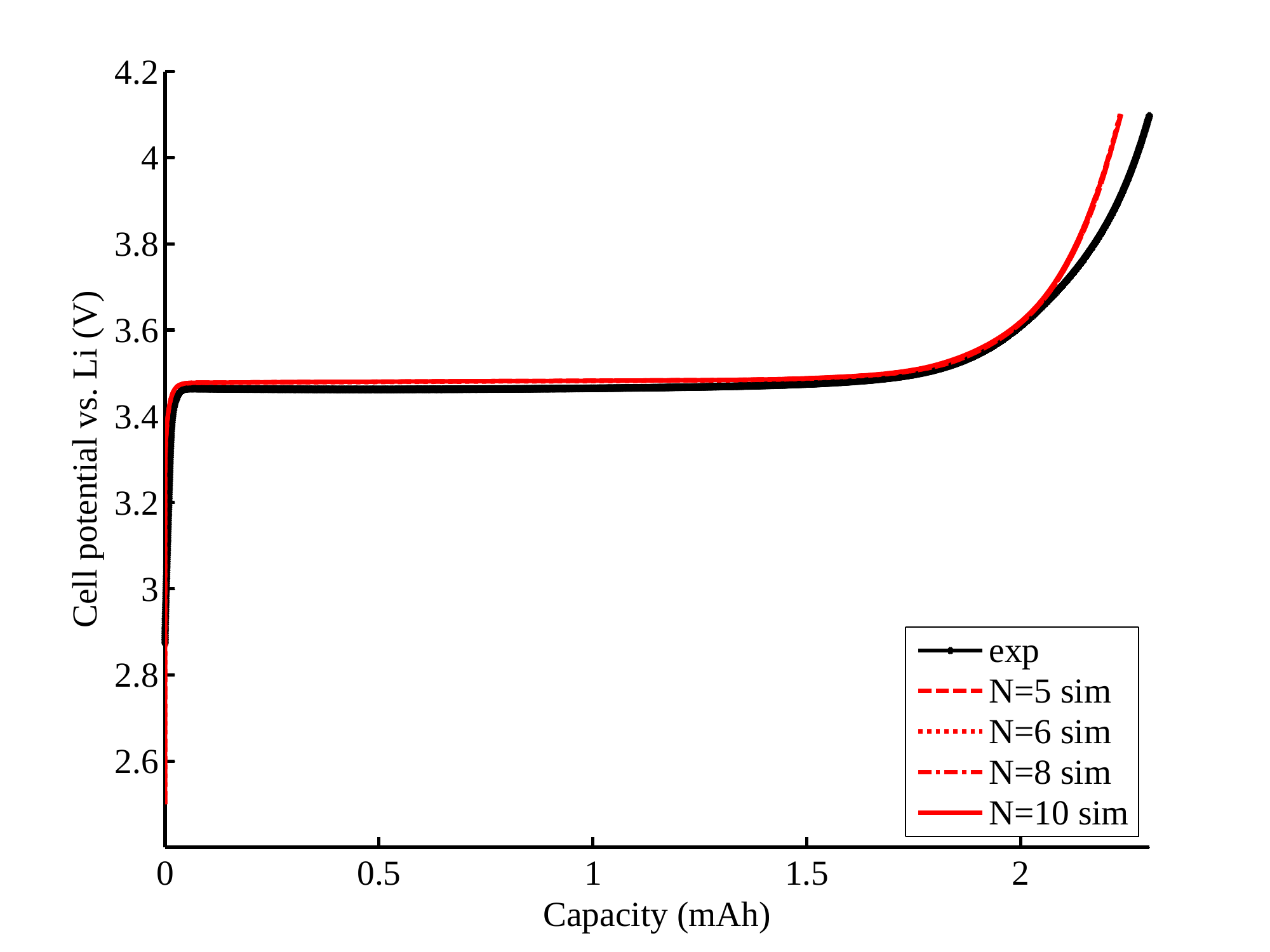}
}
\subfloat[Charging current rate=0.5C]{
  \includegraphics[width=65mm]{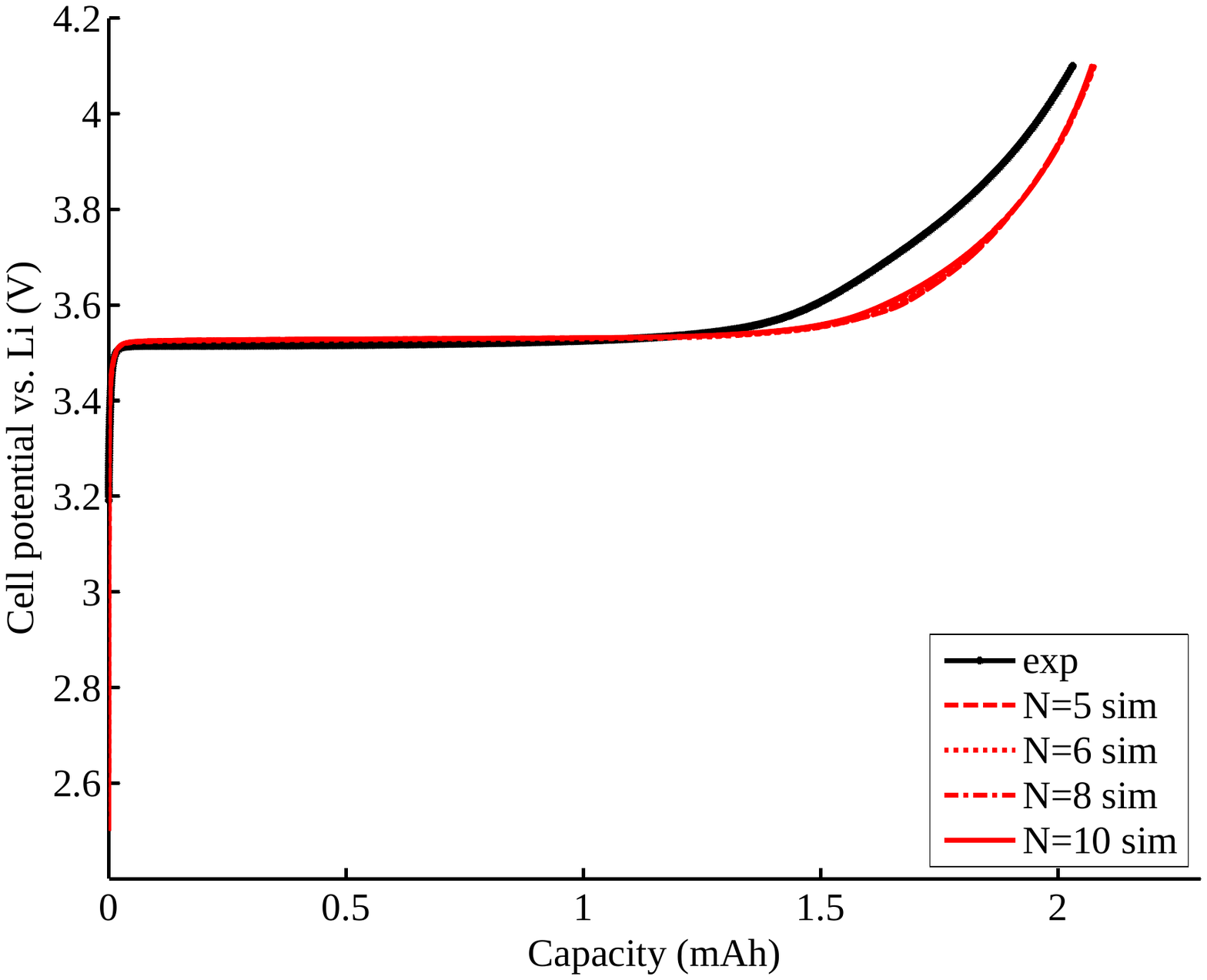}
}
\hspace{0mm}
\subfloat[Charging current rate=1C]{
  \includegraphics[width=65mm]{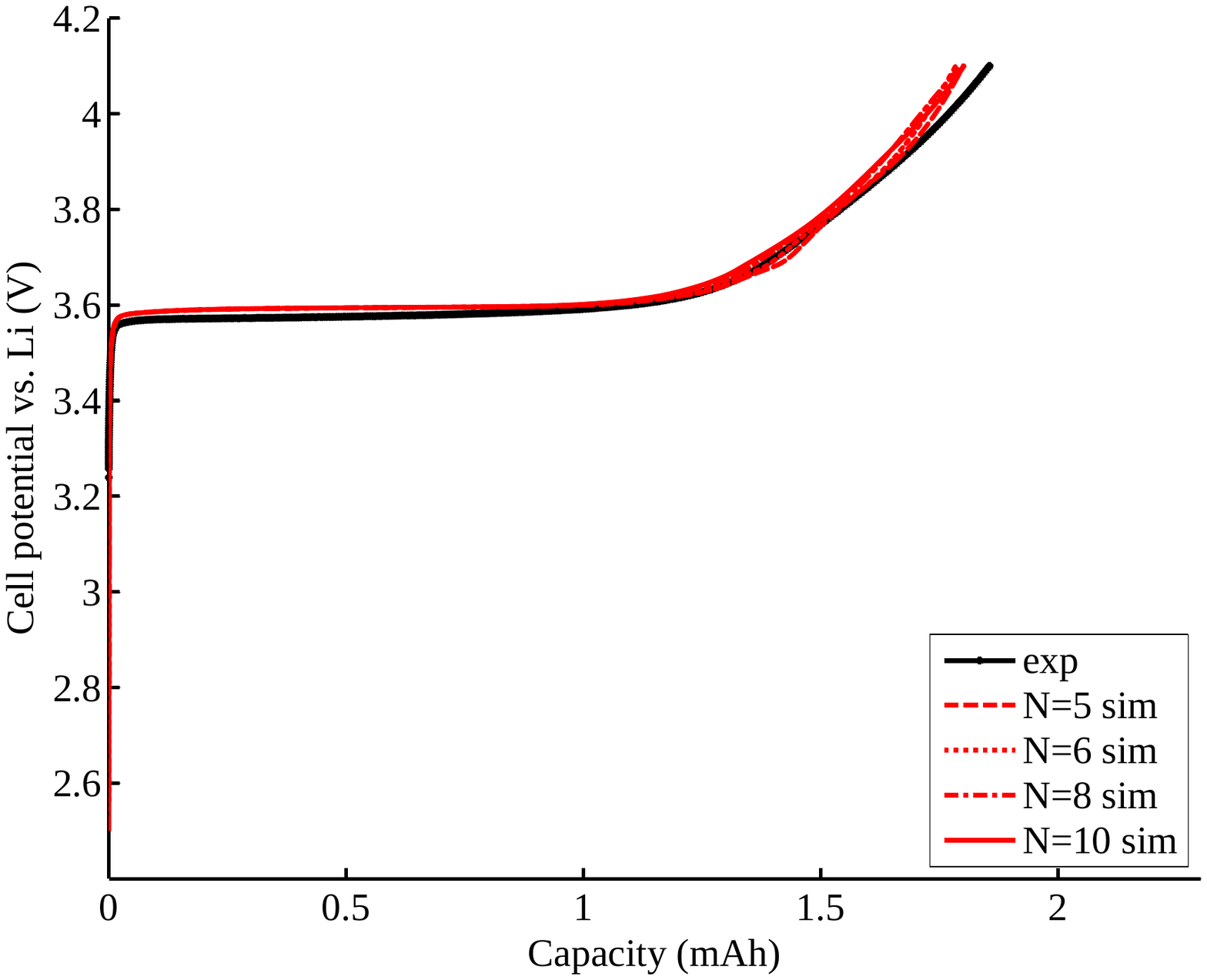}
}
\subfloat[Discharging current rate=0.2C]{
  \includegraphics[width=65mm]{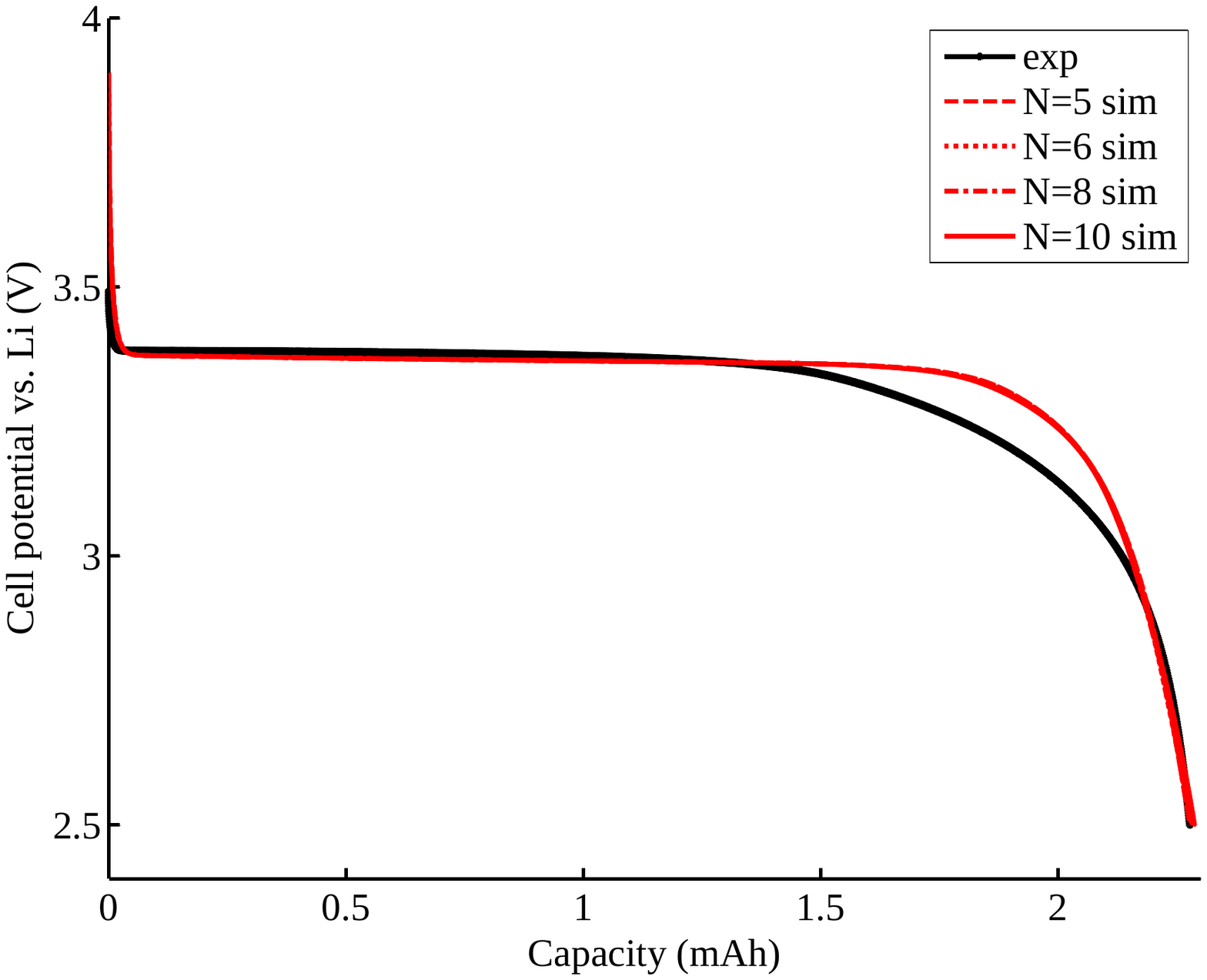}
}
\hspace{0mm}
\subfloat[Discharging current rate=0.5C]{   
  \includegraphics[width=65mm]{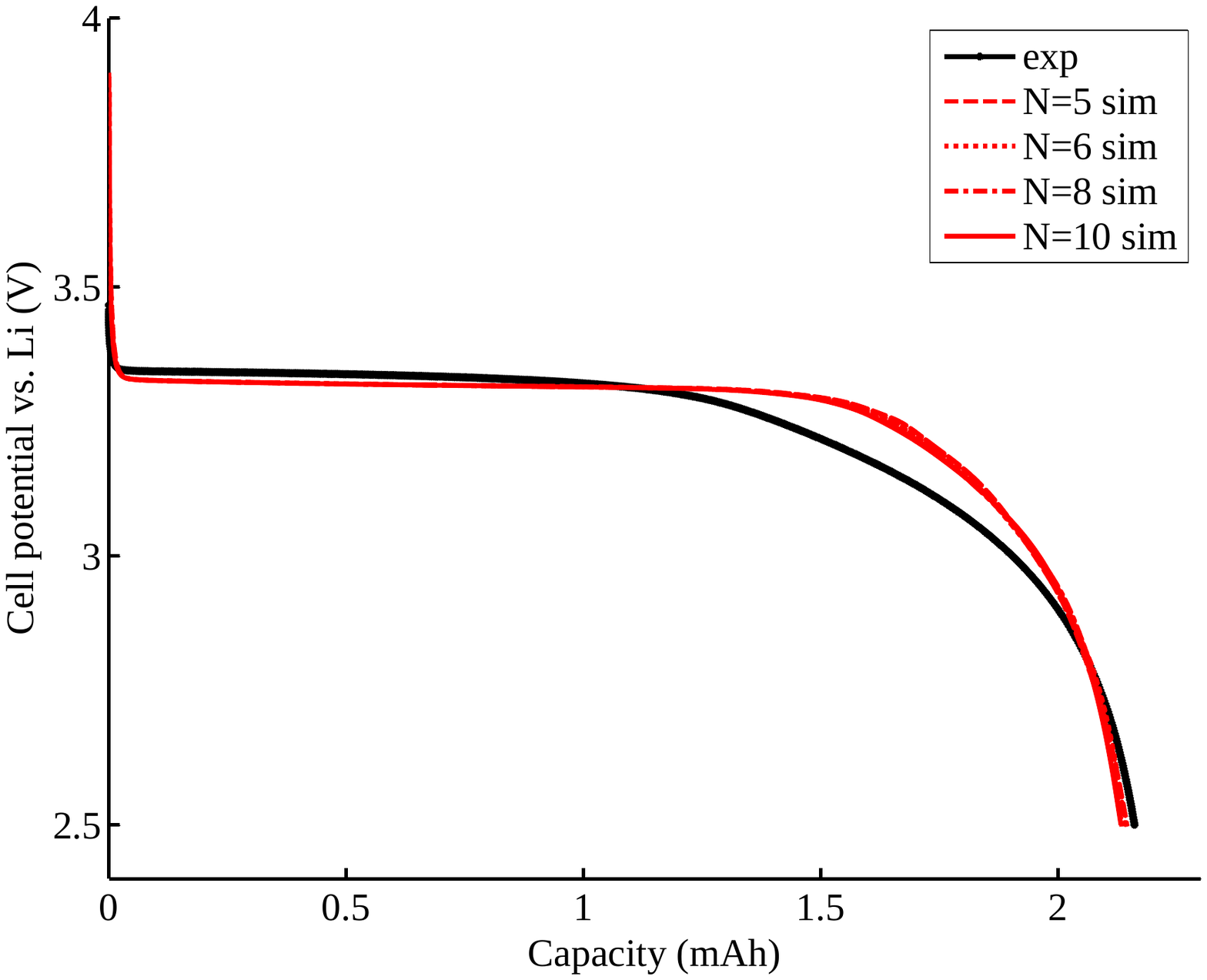}
}
\subfloat[Discharging current rate=1C]{
  \includegraphics[width=65mm]{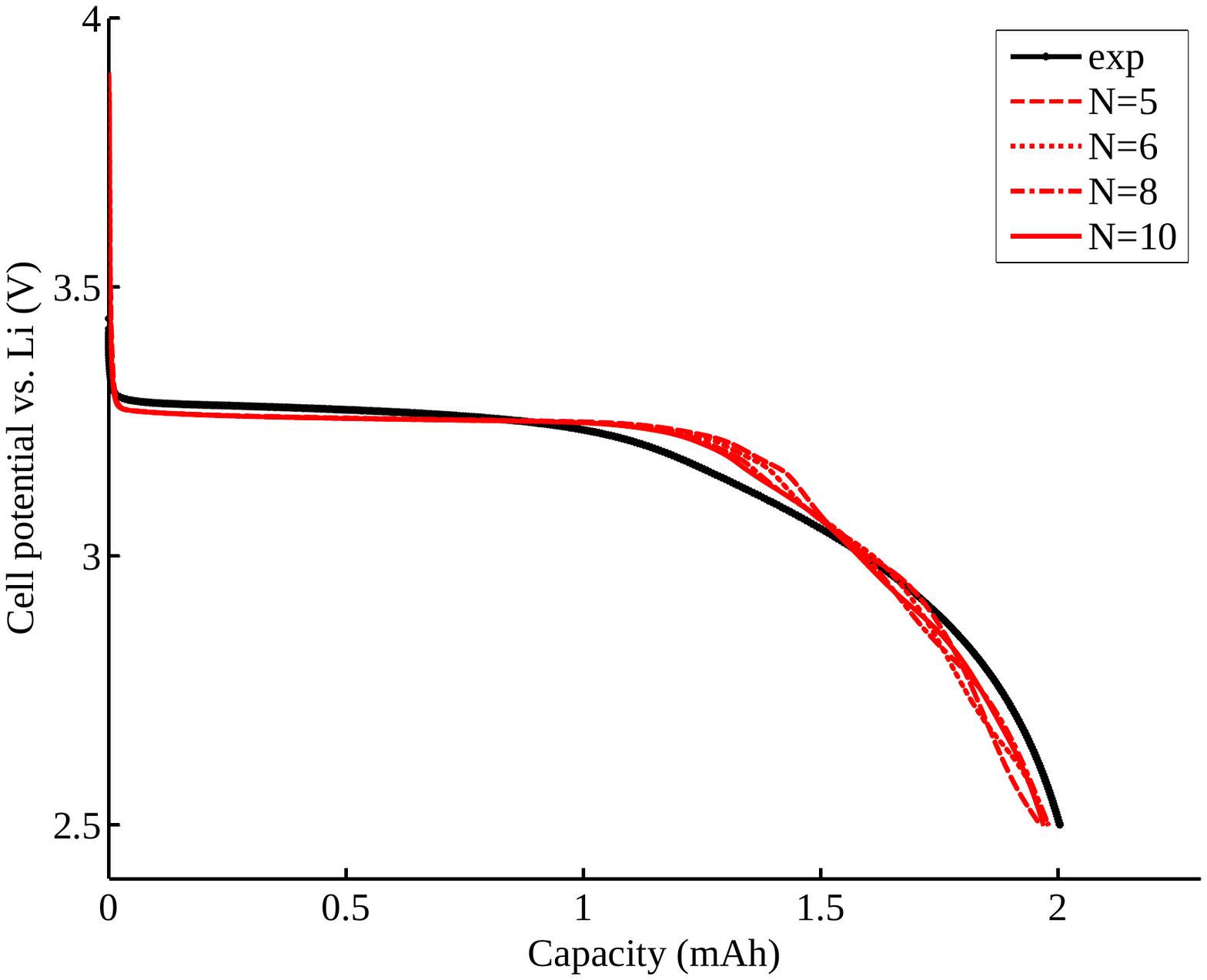}
}
\caption{Approximate output solution to the electrochemical equations at different charging/discharging current rates with the correction sampling time $Dt=3\,\text{s}$. The solutions converge for a low-order of approximation: $N\leq 8$ for the current rate $1\text{C}$ and $N\leq 6$ for the current rates less than $1\text{C}$.}
\label{F7}
\end{figure}

\end{document}